\documentclass[10pt]{amsart}
\usepackage{amsfonts,amscd}

\usepackage[arrow,matrix,curve,cmtip,ps]{xy}
\newtheorem{theorem}{Theorem}[section]
\newtheorem{lemma}[theorem]{Lemma}
\newtheorem{corollary}[theorem]{Corollary}
\newtheorem{proposition}[theorem]{Proposition}
\theoremstyle{remark}
\newtheorem{remark}[theorem]{Remark}
\theoremstyle{definition}
\newtheorem{definition}[theorem]{Definition}

\numberwithin{equation}{section} \makeatother

\begin{document}

\title{A Morita theorem for dual operator algebras}

\author{Upasana Kashyap}
\address{Department of Mathematics, University of Houston,
Houston, TX  77204-3008}
\email{upasana@math.uh.edu}
\date{\today}

\keywords{$W^{*}$-algebra; Operator algebra; Dual operator algebra;
 Dual operator module; Morita equivalence}

\begin{abstract}
We prove that two dual operator algebras are weak$^*$ Morita equivalent
in the sense of \cite{BK1}
if and only if they have equivalent categories of dual operator modules 
via completely contractive functors which are also weak$^*$-continuous on appropriate morphism spaces. 
Moreover, in a fashion  similar to the operator algebra case,
we characterize such functors as the module normal Haagerup tensor product with an appropriate weak$^*$ Morita equivalence bimodule. 
We also develop the theory of the $W^*$-dilation, which connects
the non-selfadjoint dual operator algebra with the $W^*$-algebraic framework.
In the case
of weak$^*$ Morita equivalence, 
this $W^*$-dilation is a $W^*$-module over a von Neumann algebra
generated by the non-selfadjoint dual operator algebra.
The theory of the $W^*$-dilation is a key part of the proof of our main theorem.
\end{abstract}

\maketitle
\section{Introduction and Notation}
An important and well-known perspective of understanding an algebraic object is 
to study its category of representations. For example, modules 
correspond to  
 representations of a ring hence rings are commonly  studied in terms of their modules. 
Once we view an algebraic object in terms of its category
of representations, it is natural to compare such categories.
This leads to the notion of Morita equivalence. 
The notion of Morita equivalence of rings arose in
pure algebra around 1960.
Two rings 
are defined to be Morita equivalent if and only they have equivalent categories of modules.
Morita equivalence is a powerful tool
in pure algebra, and it has inspired similar notions in operator algebra theory. In 
the 1970's Rieffel
introduced and developed the notion of Morita equivalence for $C^*$-algebras
and $W^*$-algebras.
This is a useful and  very important tool in modern operator theory.
With the advent of operator space theory in the 1990's, Blecher, Muhly and Paulsen generalized Rieffel's $C^*$-algebraic notion
of Morita equivalence to non-selfadjoint operator algebras. 


Recently we generalized Rieffel's variant of $W^*$-algebraic Morita equivalence
to dual operator algebras. 
By a {\em dual operator algebra}, we mean
a unital weak$^*$-closed algebra of 
operators on a Hilbert space which is
not necessarily selfadjoint.
One can view a dual operator algebra as a 
non-selfadjoint analogue
of a von Neumann algebra. By a non-selfadjoint
version of Sakai's theorem (see e.g.\ Section ~ 2.7 in \cite{DBbook}),
a dual operator algebra is characterized 
as a unital operator algebra that is also a
dual operator space. 

In \cite{BK1} we defined two dual
operator algebras  $M$ and $N$ to be
weak$^*$ Morita equivalent
if there exists a  dual operator $M$-$N$-bimodule $X$, 
and a dual operator $N$-$M$-bimodule $Y$,
such that $M \cong X \otimes^{\sigma h}_{N} Y$ as dual operator $M$-bimodules
(that is, via a completely isometric, weak$^*$-homeomorphism
which is also a $M$-bimodule map),
and $N  \cong Y \otimes^{\sigma h}_{M} X$  as
dual operator $N$-bimodules. 
Another notion of Morita equivalence for dual operator algebras
introduced by Eleftherakis will be briefly discussed 
later in the introduction.

In the literature of Morita equivalence for rings in pure algebra, there is a popular collection of theorems
known as  Morita I, II and III. 
Morita I can be described as the consequences of a pair of bimodules being
mutual inverses $(X \otimes_{N} Y \cong M$ and $Y \otimes_{M} X \cong N)$.
For dual operator algebras,
most of the appropriate version of Morita I is proved in \cite{BK1}.
Morita II characterizes module category equivalences 
as tensoring with an invertible bimodule, and our
main theorem here is a Morita II theorem for dual operator algebras.
The Morita III  
theorem states that there is a bijection between the set 
of isomorphism classes of invertible bimodules and the set of equivalence
classes of category equivalences;  its appropriate version
for dual operator
algebras follows as in pure algebra and will be presented in
\cite{UK}.

In \cite{BK1} we proved that two dual operator
algebras that are weak$^*$ Morita
equivalent in our sense have equivalent categories of dual operator modules.
In the present work, we prove the converse, a Morita II theorem: if two
dual operator algebras have equivalent categories of dual operator modules
then they are weak$^*$-Morita equivalent in the sense of \cite{BK1}. 
The functors implementing the categorical equivalences are characterized 
as the module normal Haagerup tensor product with an appropriate weak$^*$ Morita equivalence bimodule. 
 In Section 2, we develop the theory of the $W^*$-dilation, which connects
the non-selfadjoint dual operator algebra with the $W^*$-algebraic framework.
In particular, we use the  maximal $W^*$-algebra $\mathcal{C}$ generated by a
dual operator algebra $M$. Every dual operator $M$-module
dilates to a dual operator module over $\mathcal{C}$ which is called
the `{\em maximal dilation}'. We show that every dual operator module
is a weak$^*$-closed submodule of its maximal dilation.
Indeed, in the case
of weak$^*$ Morita equivalence this maximal dilation turns out to be a $W^*$-module ove $\mathcal{C}$.
The theory of the $W^*$-dilation is a key part of the proof of our main theorem.
In Section 3, we discuss some 
weak$^*$ Morita equivalence and $W^*$-dilation results.
In Section 4 and Section 5, we prove our main theorem.

Many
of the techniques and ideas in this paper
are taken from \cite{Max}, \cite{DB4}, \cite{DB3}, \cite{BMN}. 
We refer the reader to these papers for earlier ideas, proof techniques,
and additional details.
In some places
we just need to modify the arguments in the 
present setting of weak$^*$-topology, or merely change the tensor product.
However, 
we need to develop  new techniques to
deal with a number of subtleties that arise in
the weak$^*$-topology setting.


Another notion of Morita
equivalence for dual operator algebras was considered in \cite{Elef}  and is called $\Delta$-equivalence.
In \cite{EP} it was shown that the $\Delta$-equivalence
implies weak$^*$ Morita equivalence in the sense of \cite{BK1}. 
That is, any of the equivalences of \cite{Elef} is one of our weak$^*$ Morita equivalence.
Both the theories have different advantages.
For example, the equivalence considered in \cite{Elef} is equivalent to the
very important notion of weak$^*$ stable isomorphism.
On the other hand, our theory contains all examples considered up to this point 
in the literature of Morita-like equivalence in a dual
(weak$^*$) setting.
There are certain important examples 
that do not seems to be contained the other theory 
but are weak$^*$ Morita equivalent in our sense.
For example, in the selfadjoint setting the second dual of
strongly Morita equivalent $C^*$-algebras
are Morita equivalent in Rieffel's $W^*$-algebraic sense.
In the non-selfadjoint case, the second dual of
strongly Morita equivalent operator algebras in
the sense of Blecher, Muhly and Paulsen
are weak$^*$ Morita equivalent in our sense.
Also, two `similar' separably acting nest algebras are 
Morita equivalent in our sense but are not $\Delta$-equivalent.

In \cite{BK1} we showed that weak$^*$ Morita equivalent dual operator algebra
have equivalent categories of normal Hilbert space representations 
(also known as normal Hilbert modules). However,
the converse of this is still an open problem and at present we are working on this aspect.
The characterization theorem
in \cite{Elef} is 
in terms of equivalence of categories of normal Hilbert modules,
which intertwines  not only the representations of the
dual operator algebras, but also their restrictions to the diagonals.

We  assume that the reader is familiar with the notions 
from {\em operator space theory}. One can refer to
\cite{DBbook} and \cite{ER1} for background and most of the terminology used in this paper. We also assume that 
all dual operator algebras are unital; that is, they each have
an identity of norm 1.
We will often abbreviate `weak$^*$' to `$w^*$'.
 We reserve
the symbols $M$ and $N$ for dual operator algebras.
A {\em normal representation} of $M$ is a $w^*$-continuous
unital completely contractive homomorphism
$\pi : M \to B(H)$.
For a
dual space $X$, we let $X_*$ denote its predual.
We assume that the reader is familiar with
the weak$^*$-topology and basic duality principles
such as the Krein-Smulian theorem  (see Theorem A.2.5 in \cite{DBbook}).

A {\em concrete dual operator $M$-$N$-bimodule} is a $w^{*}$-closed
 subspace $X$ of $B(K,H)$ such that $\theta(M) X \pi(N)$ $\subset X$, where
 $\theta$ and $\pi$ are normal
 representations of $M$ and $N$ on $H$ and $K$ respectively. An
 {\em abstract
 dual operator $M$-$N$-bimodule} is defined to be an
  operator
 $M$-$N$-bimodule $X$ 
(by which we mean that $X$ is an operator space and a 
nondegenerate $M$-$N$-bimodule such that the module 
actions are completely contractive in the sense of 3.1.3 in \cite{DBbook}), 
which is also a dual operator space, such that
 the module actions are separately weak* continuous.
Such spaces can be represented completely isometrically 
as concrete dual operator bimodules (see e.g. Theorem 3.8.3 in \cite{DBbook},
\cite{ER2})
We shall write $_{M}
 \mathcal{R}$ for the category of left dual operator modules over
 $M$. The morphisms in $_{M} \mathcal{R}$ are the $w^*$-continuous
 completely bounded $M$-module maps.

By $_{M} \mathcal{H}$, we mean the category of completely
contractive normal Hilbert modules over a dual operator algebra $M$.
That is, elements of $_{M} \mathcal{H}$ are pairs $(H, \pi)$, where
$H$ is a (column) Hilbert space (see e.g.\ 1.2.23 in
\cite{DBbook}), and $\pi : M \to B(H)$ is a
 normal representation of $M$.
The module action is expressed through the equation $m \cdot \zeta =
\pi(m) \zeta$. The morphisms are bounded linear transformations
between Hilbert spaces that intertwine the representations; i.e.,\ if
$(H_{i}, \pi_{i})$, $i=1,2$, are objects of the category  $_{M}
\mathcal{H}$, then the space of morphisms is defined as:

$B_{M}(H_{1},H_{2}) = \{T \in B(H_{1},H_{2}) : T\pi_{1}(m)=
\pi_{2}(m) T $ for all $m \in M\}$. 

\noindent Any  $H \in \, _{M} \mathcal{H}$ (with its column Hilbert space
structure) is a left dual operator $M$-module.
If $E$ and $F$ are sets, then $EF$ denotes the norm
closure of the span of products $xy$ for $x \in E$ and $y \in F$.

If $X$ and $Y$ are dual operator spaces, we denote
by $CB^{\sigma}(X,Y)$ the space of completely bounded 
$w^*$-continuous linear maps from $X$ to $Y$.
Similarly if $X$ and $Y$ are left dual operator $M$-modules,
then  $CB^{\sigma}_{M}(X, Y)$ denotes the space
of completely bounded $w^*$-continuous left $M$-module
maps from $X$ to $Y$. 

 If $M$ is a  dual operator algebra, 
 then a {\em $W^*$-cover} of $M$ is a pair
 $(A,j)$ consisting of a $W^*$-algebra
 $A$ and a completely isometric $w^*$-continuous homomorphism
 $j : M \to A$, such that $j(M)$ generates $A$ as a $W^*$-algebra. 
 By the Krein-Smulian theorem $j(M)$ is a $w^*$-closed subalgebra of $A$.
 The {\em maximal
$W^*$-cover} $W^*_{\rm max}(M)$ is a $W^*$-algebra containing $M$ as
a $w^*$-closed subalgebra that is generated by $M$
as a $W^*$-algebra, and 
has the following universal property: any normal
representation $\pi : M \rightarrow B(H)$
extends uniquely to a (unital) normal $*$-representation
$\tilde{\pi} : W^*_{\rm max}(M) \rightarrow B(H)$ (see  \cite{BSo}).

 We will refer
to Rieffel's $W^*$-algebraic Morita equivalence (see \cite{R}) as
`{\em weak Morita equivalence}' for $W^*$-algebras, and the associated equivalence bimodules
as `{\em $W^*$-equivalence-bimodules}'  (see e.g. Section 8.5 in \cite{DBbook}). 

We use the normal module Haagerup 
tensor product $\otimes_{M}^{\sigma h}$ throughout the paper.
We refer to \cite{EP} and
[4, Section 2]  for the universal property 
and  general facts and properties
of  $\otimes_{M}^{\sigma h}$. Loosely speaking, the normal
module Haagerup tensor product linearizes
completely contractive balanced separately 
weak$^*$-continuous bilinear maps.
 
This work will constitue a part of the authors Ph.D. thesis
at the University of Houston. 

\section{Dual operator modules over a generated $W^{*}$-algebra and $W^{*}$-dilations}

        We begin this section with a weak$^{*}$-topology
version of Theorem 3.1 in  \cite{Max}.

\begin{theorem} \label{one*}
Let $D$ be a $W^{*}$-algebra, $B$ a Banach algebra which is 
also a dual Banach space, and $\theta : D \to B$ a unital 
$w^{*}$-continuous contractive homomorphism. 
Then the range of $\theta$ is  $w^{*}$-closed, 
and possesses an involution with respect to 
which $\theta$ is a $*$-homomorphism and 
the range of $\theta$ is a $W^{*}$-algebra.
\end{theorem}
 
 \begin{proof}
 It is known that  (see e.g.\ Theorem A.5.9 in \cite{DBbook})
the range
of a  contractive
homomorphism between a $C^{*}$-algebra
and a Banach algebra is a $C^{*}$-algebra and 
moreover such homomorphisms are $*$-homomorphisms. 
To see that the range of $\theta$ is $w^{*}$-closed, consider the quotient
 map $D/ker(\theta) \to B$  which is an isometry, and  apply the Krein-Smulian theorem.
 \end{proof}
 
 Thus if $X$ is a left dual operator module over a 
 $W^{*}$-algebra $D$, and if we let $\theta : D \to CB(X)$
 be the associated unital $w^{*}$-continuous  
contractive (equivalently
completely 
contractive by Proposition 1.2.4 in \cite{DBbook}) homomorphism, 
 then the range of $\theta$ is a $W^{*}$-algebra.
 
 \begin{theorem} \label{two*}
 Suppose that $X$ is a dual operator module over a dual 
 operator algebra $M$. Let $\theta : M \to CB(X)$ be the 
associated completely
 contractive  homomorphism. Suppose that $D$ is any $W^{*}$-algebra 
 generated by $M$. Then the $M$-action on $X$ can be 
 extended to a $D$-action with respect to which $X$ is a dual
 operator $D$-module if and only if $\theta$ is the restriction
 to $M$ of a $w^{*}$-continuous contractive 
(equivalently completely contractive)
 homomorphism 
 $\phi : D \to CB(X)$. This extended $D$-action, or equivalently
 the homomorphism $\phi$, is unique if it exists.       
 \end{theorem} 

\begin{proof}
If $\theta$ is the restriction
 to $M$ of a $w^{*}$-continuous completely contractive homomorphism 
 $\phi : D \to CB(X)$ then 
 the $M$-action on $X$ can be 
 extended to a $D$-action via $d \cdot x = \phi(d) \cdot x$. 
Note that the $D$-module action $x \mapsto d  x$ on $X$,
for $x \in X$ and $d \in D$, is a multiplier (see e.g.\ Section 4.5 in \cite{DBbook}), hence 
it is weak$^*$-continuous by Theorem 4.1 in \cite{BM2}.
The $D$-module action on $X$ is separately $w^*$-continuous
and completely contractive. Hence $X$ is a dual operator $D$-module.
The converse is obvious.
To see the uniqueness assertion, suppose that
$\phi_{1}$ and $\phi_{2}$ are two $w^{*}$-continuous 
contractive homomorphisms  $D \to CB(X)$, extending $\theta$.
By Theorem \ref{one*}, the ranges $\mathcal{E}_{1}$ and 
$\mathcal{E}_{2}$, of $\phi_{1}$ and $\phi_{2}$ respectively, are each $W^{*}$-algebras, 
but with possibly different involutions and weak$^{*}$-topologies.  
We will write these involutions
as $\star$ and $\#$ respectively.
With respect to these involutions
$\phi_{1}$ and $\phi_{2}$ are $*$-homomorphisms.
Note, $CB(X)$ is a unital Banach algebra and $\mathcal{E}_{1}$
and $\mathcal{E}_{2}$ may be viewed as unital subalgebras 
of $CB(X)$, with the same unit.
Let $a$ $\in$ $M$ and $f$ be a state on $CB(X)$. Then $f|{\mathcal{E}_{i}}$
is a state on $\mathcal{E}_{i}$ for $i = 1, 2$.
Thus $f(\phi_{1}(a)^\star) = \overline{f(\phi_{1}(a))} = \overline {f(\phi_{2}(a))} = f(\phi_{2}(a)^{\#} )$.
Thus $u = \phi_{1}(a)^\star - \phi_{2}(a)^{\#}$ is a Hermitian element
in $CB(X)$ with numerical radius 0, hence $u$ = $0$.
This implies that $\phi_{1}(a^*)  = \phi_2(a^{*})$, 
since $\phi_1$ and $\phi_2$ are $*$-homomorphisms. 
Hence $\phi_1$ equals $\phi_2$ on the $*$-subalgebra generated by $M$ in $D$.
By weak$^{*}$-density, it follows that $\phi_1$ = $\phi_2$ on $D$.
 \end{proof}
This immediately gives the following:

\begin{corollary} \label{cor}
Let $D$ be a $W^{*}$-algebra generated by a dual operator 
algebra $M$. If $X_{1}$ and $X_{2}$ are two dual operator 
$D$-modules, and if \  $T : X_{1} \to X_{2}$ is a 
$w^{*}$-continuous completely isometric and surjective $M$-module map,
then $T$ is a $D$-module map.
\end{corollary}

\begin{corollary}
Let $D$ be a $W^{*}$-algebra generated by a dual operator 
algebra $M$. Then the category $_{D} \mathcal{R}$ 
of dual operator modules over $D$
is a subcategory of the category  $_{M} \mathcal{R}$ 
of dual operator modules over
$M$. Similarly, $_{D} \mathcal{H}$ is a subcategory of $_{M} \mathcal{H}$.
\end{corollary}

Next we discuss 
the  $W^*$-dilation which we call the {\em `D-dilation'} 
of a dual operator $M$-module $X$, 
where $D$
is a $W^{*}$-algebra generated by $M$. Strictly speaking, it should be called $W^*$-$D$-dilation,
 but for brevity we will 
use the shorter term.

\begin{definition} \label{dilation}
A pair $(E,i)$ is said to be a {\em $D$-dilation} 
of a left dual operator $M$-module $X$, if the following hold:
\begin{enumerate}
\item $E$ is a left dual operator $D$-module 
and $i : X \to E$ is 
a $w^{*}$-continuous completely contractive $M$-module map. 

\item For any left dual operator $D$-module $X'$, and any 
$w^{*}$-continuous completely bounded $M$-module map 
$T : X \to X'$, there exists a unique $w^{*}$-continuous completely 
bounded $D$-module map 
$\tilde{T} : E \to X'$ such that $\tilde{T} \circ i$
= $T$, and also 
$\lVert T\rVert_{cb}$ = $\lVert \tilde{T} \rVert_{cb}$. 
\end{enumerate}
Some authors also use the terminology `{\em $D$-adjunct}' for $D$-dilation (see \cite{Max}).
\end{definition}  

The assertion in $(2)$ above implies 
that $i(X)$ generates $E$ as a dual
operator $D$-module. To see this, let $E'$ = $\overline{Di(X)}^{w*}$, 
and consider the quotient map $ q : E \to E/E'$. Then $E/E'$ is
a left dual operator $D$-module such that $q \circ i$ = $0$. Hence the assertion  in
$(2)$ in the above definition implies 
that the map $q$ = 0. Thus $E$ = $E'$.

Up to a complete isometric module isomorphism
there is a unique pair $(E,i)$ satisfying $(1)$ and $(2)$
in the above definition. To see this, 
let $(E',i')$ be any other
pair satisfying $(1)$ and $(2)$, then there 
exists a unique $w^{*}$-continuous 
completely contractive $D$-module linear maps 
$\rho : E \to E'$  and $\phi : E' \to E$
such that $\rho \circ i$ = $i'$ and $\phi \circ i'$ = $i$.
One concludes that $\rho \circ \phi$ is the
identity map on $i'(X)$ and  $\phi \circ \rho$
is the identity map on $i(X)$. Since $i(X)$ and $i'(X)$
generate $E$ as a dual operator $D$-module, and since $\phi$ and $\rho$
are $w^{*}$-continuous complete contractions, this implies that 
$\phi$ and $\rho$ are 
complete isometries. 

\begin{remark}
From the above it is clear that 
the $D$-dilation $(E,i)$ is the
unique pair satisfying $(1)$, and
such that
for all dual operator $D$-modules
$X'$, the canonical map 
$i^{*} :  CB^{\sigma} _{D}(E,X') \to   CB^{\sigma} _{M}(X,X')$,
given by composition with $i$, is an  isometric isomorphism. 
Note that 
by using $(1.7)$ and Corollary 1.6.3 in \cite{DBbook},
it is easy to see that
$M_{n}( CB^{\sigma}(X, Y))$ $\cong$ $  CB^{\sigma}(X, M_n(Y))$ 
completely isometrically for dual operator spaces $X$ and $Y$.
If $X$ is a left dual operator $M$-module, then
$M_{n}(X)$ is also a left dual operator $M$-module
via $m \cdot [x_{ij}] = [m \cdot x_{ij}] = I_{n} \otimes m  \cdot  [x_{ij}]$,
where $ I_{n} \otimes m$ denotes the diagonal matrix in $M_{n}(M)$
with diagonal entries $m$.
Indeed, if $X$ is a dual operator $M$-module,
the above module action is completely contractive 
and by Corollary 1.6.3 in \cite{DBbook}, this action is separately
$w^*$-continuous. This proves that $M_n(X)$ is  
a dual operator $M$-module if $X$ is a dual operator $M$-module. 
Since $i^{*}$ is an isometry for all dual operator $D$-modules $X'$,
it follows that $  CB^{\sigma}_{D}(E, M_n(X'))$ $\cong   CB^{\sigma}_{M}(X, M_n(X'))$ 
for all  dual operator $D$-modules $X'$,
which implies that $i^*$ is a complete isometry.
Thus the $D$-dilation $E$ of $X$ satisfies:
\begin{equation} \label{dagger-eqn}
 CB^{\sigma} _{D}(E,X') \  \cong  \      CB^{\sigma} _{M}(X,X')   
\end{equation} 
\noindent completely isometrically.
\end{remark}
\medskip
By the dual operator module version of Christensen-Effros-Sinclair  
theorem (see e.g.\ Theorem 3.3.1 in \cite{DBbook}), 
 $X'$ in Definition \ref{dilation}
can be taken to be $B(H,K)$, where $K$ is a normal
Hilbert $D$-module and $H$ is a Hilbert space.
In fact, by a modification of Theorem 3.8 in \cite{Max}, we may 
take $X' = K$. We are going to prove this
important fact in the next theorem but before
that we need to recall some tensor products facts.

For operator spaces
$X$ and $Y$, we denote the Haagerup tensor product of
$X$ and $Y$ by $X \otimes_{h} Y$. If $Z$ is another
operator space, $CB(X \times Y ,Z)$
denotes the space of completely bounded bilinear maps 
from $X \times Y \to Z$
(in the sense of Christensen and Sinclair).
It is well known that $CB(X \times Y, Z)$ $\cong $ $CB( X \otimes_{h} Y, Z) $
completely isometrically (see e.g. 1.5.4 in \cite{DBbook}).

If $X$ and $Y$ are two dual operator spaces, we 
use $(X \otimes_h Y)_{\sigma}^* $ to denote the 
subspace of $(X \otimes_{h} Y)^*$ corresponding to the
completely bounded bilinear maps from $ X \times Y \to \mathbb{C}$
which are separately $w^*$-continuous. Then we define
the {\em normal Haagerup tensor product} $X \otimes^{\sigma h} Y$
to be the operator space dual of $(X \otimes_h Y)_{\sigma}^* $.
If $Z$ is another dual operator space,
we denote by $CB^{\sigma}(X \times Y, Z)$  
the space of completely bounded bilinear maps
from $X \times Y \to Z$ 
which are separately $w^*$-continuous.
By the matrical version of $(5.22)$ in \cite{ER}, 
$CB^{\sigma}(X \times Y, Z)$ $\cong$ $CB^{\sigma}(X \otimes^{\sigma h} Y, Z)$
completely isometrically.

Suppose $X$ is a right dual operator $M$-module
and $Y$ is a left dual operator $M$-module. 
A bilinear map $u : X \times Y \to Z$
is {\em $M$-balanced} if $u(xm,y) = u(x,my)$ for $m \in A$.
We let 
$(X \otimes_{hM} Y)_{\sigma}^* $ denote the 
subspace of $(X \otimes_{hM} Y)^*$ corresponding to the
completely bounded balanced bilinear maps from $ X ~ \times ~ Y  ~ \to ~ \mathbb{C}$
which are separately $w^*$-continuous, where $\otimes_{hM}$
denotes the module Haagerup tensor product 
(see e.g. 3.4.2, 3.4.3 in \cite{DBbook}). By Proposition 2.1
in \cite{EP}, the {\em  module normal Haagerup tensor product}
$X \otimes^{\sigma h}_{M} Y$
may be defined to be the operator space dual of 
$(X \otimes_{hM} Y)_{\sigma}^* $.
If $Z$ is another dual operator space,
we denote by $CB^{ M \sigma }(X \times Y, Z)$  
the space of completely bounded balanced  
separately $w^*$-continuous bilinear maps.
By Proposition 2.2 in  \cite{EP}, 
$CB^{ M \sigma }(X \times Y, Z)$ $\cong$ $CB^{\sigma}(X \otimes^{\sigma h}_{M} Y, Z)$
completely isometrically.

In order to prove the next lemma, we will introduce
some notation.
Let $CB^{ S \sigma  }(X \buildrel  \frown  \over \otimes Y , Z)$ 
denote the subspace
of $CB(X \buildrel  \frown  \over \otimes Y, Z)$ 
consisting of completely bounded maps from $X \buildrel  \frown  \over \otimes Y $ to $Z$
that are induced by the jointly completely
bounded bilinear maps from $X \times Y \to Z$ which are separately $w^*$-continuous,
where $\buildrel  \frown  \over \otimes $ denotes the operator space projective
tensor product (see e.g. 1.5.11 in \cite{DBbook}). In the case,
when $Z = \mathbb{C}$, we denote $CB^{S \sigma  }(X \buildrel  \frown  \over \otimes Y , \mathbb{C})$ 
by $(X \buildrel \frown \over \otimes Y)_{\sigma}^{*}.$ 

\begin{lemma}
For any Hilbert spaces $H$ and $K$ and dual operator space $X$,
$ CB^{\sigma}(X, B(H,K))$ $\cong$ $ CB^{\sigma}(X \otimes^{\sigma h} H^c, K^c)$ 
$\cong$ $(\overline{K}^r \otimes^{\sigma h}  X \otimes^{\sigma h} H^c)_*$
completely isometrically. 
\end{lemma}

\begin{proof}
For any dual operator space $X$, we have the following isometries:
\begin{eqnarray*}
CB^{\sigma} (X \otimes^{\sigma h} H^c, K^c) &\cong& CB^{\sigma} (X \times H^c, K^c) \\
&\cong& CB^{S \sigma } (X \buildrel  \frown  \over \otimes H^c, K^c) \\
&\cong& CB^{\sigma} (X, CB(H^c, K^c) \\
&\cong& CB^{\sigma} ( X , B(H, K))
\end{eqnarray*}
using  Proposition 1.5.14 (1) and $(1.50)$ from \cite{DBbook}.
Consider 
\begin{eqnarray*}
 CB^{\sigma}(X \otimes^{\sigma h} H^c, K^c) 
&\cong&  (\overline{K}^r \buildrel  \frown \over \otimes  (X \otimes^{\sigma h} H^c ))_{\sigma}^*\\
&\cong&  (\overline{K}^r \otimes_{h}  (X \otimes^{\sigma h} H^c ))_{\sigma}^*\\
&\cong&  (\overline{K}^r \otimes^{\sigma h}  (X \otimes^{\sigma h} H^c))_*\\ 
&\cong&  (\overline{K}^r \otimes^{\sigma h}  X \otimes^{\sigma h} H^c)_*,
\end{eqnarray*}
using  $(1.51)$ and Proposition
1.5.14 (1) in \cite{DBbook}, and associativity of the normal
Haagerup tensor product.
\end{proof}

Similarly we have the module version of the above lemma:

\begin{lemma} \label{mm}
Let $X$ be a left dual operator $M$-module and $K$ be a 
normal Hilbert $M$-module. Then for any Hilbert space  $H$, 
$ CB^{\sigma}_{M}(X, B(H,K))$ $\cong$ $ CB^{\sigma}_{M}(X \otimes^{\sigma h} H^c, K^c)$ 
 $\cong$ $(\overline{K}^r \otimes^{\sigma h}_{M} X \otimes^{\sigma h} H^c)_* $
completely isometrically.
\end{lemma}

\begin{proof}
The first isomorphism follows as above with
completely bounded maps  replaced with module completely bounded maps.
Consider  
\begin{eqnarray*}
CB^{\sigma}_{M}(X \otimes^{\sigma h}  H^c, K^c)  
&\cong&  (\overline{K}^r \buildrel \frown \over \otimes_{M}  (X \otimes^{\sigma h} H^c ))_{\sigma}^*\\
&\cong&  (\overline{K}^r \otimes_{hM}  (X \otimes^{\sigma h} H^c ))_{\sigma}^*\\
&\cong& (\overline{K}^r \otimes^{\sigma h}_{M}  (X \otimes^{\sigma h} H^c))_* \\
&\cong& (\overline{K}^r \otimes^{\sigma h}_{M}  X \otimes^{\sigma h} H^c)_*,
\end{eqnarray*}
using  Corollary 3.5.10 in \cite{DBbook},  $K^r \otimes_{hM} - = ~ K^r  \buildrel \frown \over \otimes_{M} -$ 
and a variant of
Proposition ~ 2.9 in \cite{BK1}.
\end{proof}

We would like to thank David Blecher for the proof of the following lemma.

\begin{lemma} \label{gi}
Let  $S : X \to Y$ be a $w^*$-continuous
linear map between dual operator spaces. The following are equivalent:\\
\noindent (i) $S$ is a complete isometry and surjective.\\
\noindent (ii) For some Hilbert space $H$, $S \otimes I_H: X \otimes^{\sigma h} H^c \to
Y \otimes^{\sigma h} H^c$ is a complete isometry and surjective.  
\end{lemma}

\begin{proof}
Firstly, suppose $S$ is a completely isometric and $w^*$-homeomorphic map.
Then, by the functoriality of the normal Haagerup tensor product
$S \otimes I_H$ and $S^{-1} \otimes I_H$ are completely
contractive $w^*$-continuous maps,
 where $I_H$ denotes the identity map on $H$.
. 
Also $(S^{-1} \otimes I_H) \circ  (S \otimes I_H)
= Id$ on a weak$^*$-dense subset $X \otimes H$. 
By $w^*$-density, $(S^{-1} \otimes I_H) \circ (S \otimes I_H)
= Id$ on $X \otimes^{\sigma  h} H^c$. Similarly,
 $(S \otimes I_H)\circ (S^{-1} \otimes I_H)$
= $Id$. Thus $S \otimes I_{H}$ is a completely isometric
and $w^*$-homeomorphic map. 

Conversely, suppose $(ii)$ holds.
Fix a $\eta \in H$ with $\lVert\eta\rVert$ = 1.
Let $v : X \to X \otimes^{\sigma h} \eta$ : $x \mapsto x \otimes \eta$.
Since $ X \subseteq X \otimes_h H^c$ completely isometrically via $v$,
and $ X \otimes_h H^c \subseteq X \otimes^{\sigma h} H^c$ completely
isometrically, this implies that $v$ is a complete isometry.
If $S \otimes I_H$ is a complete isometry,
then  $S \otimes I_H$ restricted to $X \otimes^{\sigma h}  \eta$
is a complete isometry. 
Similarly, let $u : Y \to Y \otimes^{\sigma h} \eta : y \mapsto y \otimes \eta$.
Thus, $S$ =  $u^{-1} \circ ( S \otimes I_H) \circ v$
is a complete isometry. To see $S$ is onto, 
suppose for the sake of contradiction that it is not.
Then by the Krein-Smulian theorem $G$ = Ran$(S)$ is a weak$^*$-closed
proper subspace of $Y$. Let $\varphi \in G_{\perp}$  and $ \varphi \neq 0$.
Consider a map $r : Y \otimes^{\sigma h} H^c \to \mathbb{C} \otimes^{\sigma h} H^c$ : $ y \otimes \zeta \mapsto \varphi(y) \otimes \zeta$.
Then $r \circ (S \otimes I_H)$ = 0, since this vanishes  on a $w^*$-dense
subset $Y \otimes H^c$. So $r =0$. Hence $\varphi(y) \otimes \zeta = 0$
for all $\zeta \in H$ and $y \in Y$. This implies
 $\varphi = 0$, which is a contradiction.      
\end{proof}

\begin{theorem} \label{wow}
Suppose $E$ is a left dual operator $D$-module and $i : X \to E$ is 
a $w^{*}$-continuous completely contractive $M$-module map.  Then
$(E,i)$ is the $D$-dilation of $X$ if and only if the canonical map
$i^{*} :   CB^{\sigma} _{D}(E,K) \to   CB^{\sigma}_{M}(X,K)$ as  defined above
is a complete isometric isomorphism, for all normal Hilbert $D$-modules
$K$. It is sufficient to take $K$ to be the normal
universal representation of $D$ or any normal  
generator for $_{D} \mathcal{H}$
in the sense of \cite{BSo},  \cite{R}. 
\end{theorem}

\begin{proof}
Consider the following sequence of complete contractions:
\begin{center}
$\overline{K}^r \otimes^{\sigma h}_{M} X$ \  $\buildrel    {id \otimes i}  \over \longrightarrow$ \ $\overline{K}^r \otimes^{\sigma h}_{M} E$ 
$\cong$ $\overline{K}^r \otimes^{\sigma h}_{D} D \otimes^{\sigma h}_{M} E$ $ \to \overline{K}^r \otimes^{\sigma h}_{D} E$
\end{center}
where the last map in the sequence comes from the multiplication $ D \times E$ $\to E$. Taking the composition of the above maps, we get
a complete contraction $S : \overline{K}^r \otimes^{\sigma h}_{M} X \to \overline{K}^r \otimes^{\sigma h}_{ D} E$. 
Tensoring $S$ with the identity map on $H$, we get
a $w^*$-continuous,
completely contractive linear map
$S_1$ = $S \otimes id_{H} : \overline{K}^r \otimes^{\sigma h}_{M} X \otimes^{\sigma h} H^c \to \overline{K}^r \otimes^{\sigma h}_{D} E \otimes ^{\sigma h} H^c $ by Corollary
2.4 in \cite{BK1}.
From a well known weak$^*$-topology fact, $S_1$ = $T^*$
for some $T : (\overline{K}^r \otimes^{\sigma h}_{D} E \otimes ^{\sigma h} H^c)_*
\to  (\overline{K}^r \otimes^{\sigma h}_{M} X \otimes^{\sigma h} H^c)_* $. 
From Lemma ~ \ref{mm},  and standard weak$^*$-density arguments,
it follows that $T$ equals $i^*$, as defined earlier.
Indeed, we use the duality pairing, namely,
$\langle \psi \otimes x \otimes \eta, T \rangle = \langle T(x)(\eta), \psi \rangle$,
for $T \in CB^{\sigma}_{M}(X,B(H,K))$,
$x \in X$, $\eta \in H$, $\psi \in K^*$,
to check that $(i^*)^* = S_1$ on the
weak$^*$-dense subset $\overline{K}^r \otimes X \otimes H^c$.
Then by weak$^*$-density, it follows that $(i^*)^* = S_1 = T^*$, so $i^* = T$.
Hence, $i^*$ is an isometric isomorphism if and only if
$S_1$ is an isometric isomorphism  if and only if 
$S$ is an isometric isomorphism by Lemma \ref{gi}.
Note that with $H$ = $\mathbb{C}$ in Lemma \ref{mm},
$CB^{\sigma}_{M}(X, K^c)$ = $ (\overline{K}^r \otimes^{\sigma h}_{M} X)_{*} $.
From Lemma ~ \ref{mm}, it is clear
that $ CB^{\sigma}_{D}(E,K^c)$ $\cong$  $ CB^{\sigma}_{M}(X,K^c)$ if and only if
$ CB^{\sigma}_{D}(E \otimes^{\sigma h} H^c, K^c)$ $\cong$ $ CB^{\sigma}_{M}(X \otimes^{\sigma h} H^c, K^c)$ for all normal Hilbert $D$-modules $K$.
For the last assertion, note that every nondegenerate
normal Hilbert $D$-module $K$ is a complemented submodule
of a direct sum of $I$ copies of the normal universal representation
or normal generator, for some cardinal $I$ (see e.g. \cite{BSo}).
Therefore we need to show that  
if $CB^{\sigma} _{D}(E,K) \cong   CB^{\sigma}_{M}(X,K)$ completely isometrically 
then  $CB^{\sigma} _{D}(E,K^{I}) \cong   CB^{\sigma}_{M}(X,K^I )$ completely
isometrically as well, where $K^{I}$ denotes
the Hilbert space direct sum of $I$-copies of $K$.
 This follows from the operator space fact
that $CB^{\sigma} _{M}(X,Y^{I})$ $\cong$ $M_{I, 1}(CB^{\sigma} _{M}(X,Y))$
completely isometrically for any dual operator spaces $X$ and $Y$ which are
also $M$-modules (see page 156 in \cite{ER2}). Here $M_{I,1}(X)$ denotes the operator space
of columns of length $I$ with entries in $X$, whose finite subcolumns have 
uniformly bounded norm.
\end{proof}

The following lemma shows the existence
of the $D$-dilation. The normal
module Haagerup tensor product $D \otimes^{\sigma h}_{M} X$ 
(which is a dual operator $D$-module by Lemma ~ 2.3 in \cite{BK1}) acts as
the $D$-dilation of $X$. We note 
that, since by Lemma 2.10 in \cite{BK1}  $M \otimes^{\sigma h}_{M} X$ $\cong$ $X$, there is a
canonical $w^{*}$-continuous completely contractive $M$-module
map $i : X \to D \otimes^{\sigma h}_{M} X$ taking 
$x \mapsto 1\otimes_{M} x$.

\begin{lemma} \label{d}
For any left dual operator module $X$ over $M$, the dual 
operator $D$-module $E = D \otimes^{\sigma h}_{M} X$ is the
$D$-dilation of $X$.
\end{lemma}

\begin{proof}
If $T : X \to X'$ is as in Definition 
\ref{dilation}, then by the functoriality of the 
normal module Haagerup tensor product, 
$\mathbb{I}_{D} \otimes T  : D \otimes^{\sigma h}_{M} X \to  D \otimes^{\sigma h}_{M} X'$
is $w^{*}$-continuous completely bounded. Composing 
this with the $w^{*}$-continuous module action $D \otimes^{\sigma h}_{M} X' \to X'$
gives the required map $\tilde{T}$. It is routine to check that $\tilde{T}$
has the required properties. 
\end{proof}

\begin{lemma}
If $X$ is a left dual operator $M$-module, and if $D$ is a 
$W^{*}$-algebra generated by $M$, then the following 
are equivalent:
\begin{enumerate}
 \item There exists a dual operator $D$-module $X'$
 and a completely isometric $w^{*}$-continuous 
 $M$-module map $j : X \to X'$.
 
 \item The canonical $w^{*}$-continuous $M$-module
 map $i : X \to D \otimes^{\sigma h}_{M} X$, is a 
 complete isometry.
 \end{enumerate}
\end{lemma}

\begin{proof}
The one direction $(2)$ implies $(1)$ is obvious. For the difficult direction, 
suppose that $m$ is the module action on $X'$. Then we 
have the following sequence of canonical $w^{*}$-continuous
completely contractive $M$-module maps:
\begin{center}
$ X  \buildrel i \over \longrightarrow  D \otimes^{\sigma h}_{M} X \buildrel \mathbb{I} \otimes j \over \longrightarrow D \otimes^{\sigma h}_{M} X'  \buildrel {m} \over \longrightarrow X'$.
\end{center}
The composition of these maps equals $j$, which is a complete isometry. 
This forces $i$ to be a complete isometry
which proves the assertion. 
\end{proof}

In the case that $D = \mathcal{C} = W^*_{max}(M)$,
we call $ \mathcal{C} \otimes^{\sigma h}_{M} X$
the `maximal $W^*$-dilation' or `maximal dilation'. 
This is the key point in proving our main theorem (Section 4).
The reason we work mostly
with maximal dilation instead of any arbitrary dilation is 
the following result.

\begin{corollary} \label{C}
For any left dual operator $M$-module $X$, the canonical
$M$-module map $i : X \to \mathcal{C} \otimes^{\sigma h}_{M} X$
is a $w^{*}$-continuous complete isometry.
\end{corollary}

\begin{proof}
This follows from the previous result, the Christensen-Effros-Sinclair-representation 
theorem for dual
operator modules, and the fact that
every normal Hilbert $M$-module is a normal Hilbert  $\mathcal{C}$-module for 
the maximal $W^{*}$-algebra generated by $M$ (i.e. the universal property of $\mathcal{C} )$. 
\end{proof}

Hence, we may regard $X$ as a $w^{*}$-closed $M$-submodule 
of $\mathcal{C} \otimes^{\sigma h}_{M} X$.
There is a similar notion of  $W^{*}$-dilation for right dual
operator modules or dual operator bimodules. 
The results in this section carry 
through analogously to these cases. 

\section{Morita equivalence of dual operator algebras}

In this section, $M$ and $N$ are again dual
operator algebras.
We reserve the symbols $\mathcal{C}$ and $\mathcal{D}$ 
for the maximal $W^{*}$-algebras  
$W^*_{max}(M)$ and $W^*_{max}(N)$ generated by $M$
and $N$ respectively. We refer the reader to
\cite{BK1} if further background for this section is needed.
 
We begin with the following  normal Hilbert module
characterization of $W^*$-algebras which is proved
in Proposition 7.2.12 in \cite{DBbook}.
 
\begin{proposition} \label{W}
Let $M$ be a dual operator algebra. Then $M$ is a $W^*$-algebra
if and only if for every completely contractive normal 
representation $\pi : M \to B(H)$, 
the commutant $\pi(M)'$ is selfadjoint.
\end{proposition}

\begin{corollary}\label{n}
Suppose $M$ and $N$ are dual operator
algebras such that the categories  
$_{M} \mathcal{H}$ and $_{N} \mathcal{H}$ are completely isometrically
equivalent; i.e.,  there exist completely contractive functors 
$F : $ $_{M} \mathcal{H}  \to $ $_{N} \mathcal{H}$
 and $G : $ $ _{N} \mathcal{H}  \to$ $ _{M} \mathcal{H}$,
 such that $FG \cong Id$ and $GF \cong  Id$ completely isometrically,
 then: 
\begin{enumerate}
\item
If $M$ is a $W^{*}$-algebra then so is $N$.
\item
Also $_{\mathcal{C}} \mathcal{H}$ and $_{\mathcal{D}} \mathcal{H}$ 
are completely isometrically equivalent.
\end{enumerate}
\end{corollary}

\begin{proof}
Suppose $F : $ $_{M} \mathcal{H}  \to $ $_{N} \mathcal{H}$
 and $G : $ $ _{N} \mathcal{H}  \to$ $ _{M} \mathcal{H}$,
are functors as in the statement of the corollary. 
If $M$ is a $W^*$-algebra,
then for $H$ $\in$  $_{M} \mathcal{H}$,
$B_{M}(H)$ is a $W^*$-algebra by Proposition \ref{W}.
The map  $T \mapsto F(T)$ from $ B_{M}(H)$ to 
 $ B_{N}(F(H))$ is a surjective isometric homomorphism 
(see Lemma 2.2 in \cite{DB4} or Lemma \ref{B} below).
Hence by Theorem A.5.9 in \cite{DBbook},
 this is a $*$--homomorphism if $M$ is a $W^*$-algebra,
 and consequently its range $ B_{N}(F(H))$  is a $W^*$-algebra.
Thus, if $M$ is a $W^*$-algebra, then $ B_{N}(H)$ is a $W^*$-algebra
 for all normal Hilbert $N$-modules $H$. From  
 Proposition \ref{W},
it follows that $N$ is a $W^*$-algebra.
For $H$ $\in$  $_{M} \mathcal{H}$,
we have $ B_{\mathcal{C}}(H)$ is a subalgebra of $ B_{M}(H)$.
The proof that 
  $F$  restricts to a functor from
$_{\mathcal{C}} \mathcal{H}$ to $_{\mathcal{D}} \mathcal{H}$
and similar assertion for $G$,  follows identically to the $C^*$-algebra
case
(see e.g. Proposition 5.1 in \cite{Max}).
\end{proof} 


\begin{definition} \label{PP}
\begin{enumerate}
\item Suppose that $\mathcal{E}$ and $\mathcal{F}$ are weakly
Morita equivalent $W^{*}$-algebras in the sense of Rieffel \cite{R}, 
and that $Z$
is a $W^*$-equivalence $\mathcal{F}$-$\mathcal{E}$-bimodule (see 8.5.12 in \cite{DBbook}),
and that $W= \overline{Z}$ is the conjugate $\mathcal{E}$-$\mathcal{F}$-
bimodule of $Z$. Then we say that $(\mathcal{E},\mathcal{F},W,Z)$
is a {\em $W^{*}$-Morita context } (or 
{\em $W^{*}$-context} for short).\\

\item Suppose that $M$ and $N$ are dual operator
algebras, and suppose that $\mathcal{E}$ and $\mathcal{F}$ are $W^{*}$-algebras 
generated by $M$ and $N$ respectively.
Suppose that $(\mathcal{E},\mathcal{F},W,Z)$ is a $W^{*}$-Morita context,
$X$ is a $w^{*}$-closed $M$-$N$-submodule of $W$,
and $Y$ is a $w^{*}$-closed $N$-$M$-submodule
of $Z$. Suppose that the natural pairings $Z \times W \to \mathcal{F}$
and $W \times Z \to \mathcal{E}$ restrict to maps  $Y \times X \to N$,
and $X \times Y \to M$ respectively, both with $w^{*}$-dense range. Then
we say $(M,N,X,Y)$ is a {\em subcontext } of $(\mathcal{E},\mathcal{F},W,Z)$. 
If further, $\mathcal{E}$ and $\mathcal{F}$ are maximal $W^{*}$-covers
(as defined in the introduction)
of $M$ and $N$ respectively, then we say that 
$(M,N,X,Y)$ is a {\em maximal subcontext }.\\

\item A subcontext $(M,N,X,Y)$ of a $W^{*}$-Morita context 
$(\mathcal{E},\mathcal{F},W,Z)$ is {\em left 
dilatable} if $W$ is the left $\mathcal{E}$-dilation of $X$, and
$Z$ is the left $\mathcal{F}$-dilation of $Y$. In this case
we say that $M$ and $N$ are {\em left weakly subequivalent}
and $(M,N,X,Y)$ is a  {\em left subequivalence context}.

\end{enumerate}
\end{definition}
There is a similar definition and symmetric theory where
we replace the words `left' by `right' or `two-sided'.\\

\begin{remark}  
Note that (2) in the above definition implies that $X$ and $Y$
are nondegenerate dual operator modules over $M$ and $N$.
\end{remark}
\medskip

Write $\mathcal{L}^{w}$ for the set of $2\times2$ matrices
  $$\mathcal{L}^w =
\Big\lbrace \left[
\begin{array}{ccl}
a & x \\
y & b
\end{array}
\right] : \ a \in M, b\in N, x\in X, y\in Y \Big\rbrace.  $$ 
Write $\mathcal{L'}$
for the same set, but with entries from the $W^{*}$-context 
$(\mathcal{E},\mathcal{F},W,Z)$.  It is well known that $\mathcal{L'}$
is canonically a $W^{*}$- ~ algebra, called the `linking $W^{*}$-algebra'
of the $W^*$-context $(\mathcal{E},\mathcal{F},W,Z)$ (see e.g. 8.5.10 in \cite{DBbook}). Saying that $(M,N,X,Y)$
is a subcontext of $(\mathcal{E},\mathcal{F},W,Z)$ implies that $\mathcal{L}^{w}$ is a $w^{*}$-closed subalgebra of $\mathcal{L'}$. 
Thus a subcontext gives a linking dual operator algebra $\mathcal{L}^{w}$.
Clearly $\mathcal{L}^{w}$ has a unit. We shall see that $\mathcal{L}^{w}$ generates 
$\mathcal{L'}$ as a $W^{*}$-algebra.

\medskip

The proof of the following  theorem is
similar to the proof of Theorem ~ 5.2 in \cite{BK1}
with an arbitrary $W^{*}$-dilation in place of $W^{*}_{max}(M)$
and hence we omit it. 

\begin{theorem} \label{P}
Suppose that dual operator algebras $M$ and $N$ 
are linked by a $weak^{*}$Morita context 
$(M,N,X,Y)$ in the sense of \cite{BK1}.
Suppose that $M$ is represented normally and completely 
isometrically as a subalgebra of $B(H)$ nondegenerately,
for some Hilbert space $H$, and let $\mathcal{E}$ be the
$W^{*}$-algebra generated by $M$ in $B(H)$. Then
$Y \otimes^{\sigma h}_{M} \mathcal{E}$ is a right $W^{*}$-module 
over $\mathcal{E}$. Also (as in the proof of Theorem ~ 5.2 in \cite{BK1}) 
$Y \otimes^{\sigma h}_{M} \mathcal{E} \cong \overline{Y \mathcal{E}}^{w*}$ 
completely isometrically and $w^*$-homeomorphically and
hence $Y \otimes^{\sigma h}_{M} \mathcal{E}$  contains $Y$ as a $w^*$-closed
$M$-submodule completely isometrically.
Also, via this module, $\mathcal{E}$ is weakly Morita equivalent 
(in the sense of Rieffel)
to the $W^{*}$-algebra $\mathcal{F}$ generated
by the completely isometric induced normal representation of
$N$ on $Y \otimes^{\sigma h}_{M} H$.
\end{theorem} 

If $C$ is a $W^{*}$-algebra generated by $M$, then
we shall write $\mathcal{F}(C)$ for $ Y  \otimes^{\sigma h}_{M} C  \otimes^{\sigma h}_{M} X$. From an obvious modification of Theorem 5.2 in
\cite{BK1}, we have that $\mathcal{F}(C)$ is a 
$W^{*}$-algebra
 containing a copy of $N$, which is $*$-isomorphic 
and $w^{*}$-homeomorphic to $(YCX)^{-w*}$.
 The copy of $N$
may be identified with $(YMX)^{-w*}$. Thus, Theorem
\ref{P} tells us that $C$ is weakly Morita 
equivalent to $\mathcal{F}(C)$ as $W^{*}$-algebras.

Similarly, if $D$ is a $W^{*}$-algebra generated by 
$N$, then we write $\mathcal{G}(D)$ for 
$X  \otimes^{\sigma h}_{N} D  \otimes^{\sigma h}_{N} Y$.
Again $\mathcal{G}(D)$ $\cong$ $(XDY)^{-w*}$
$*$-isomorphically and $w^{*}$-homeomorphically.
By 
associativity of the module normal Haagerup tensor product  and Lemma ~ 2.10 in \cite{BK1}, $\mathcal{G}(\mathcal{F}(C)) \cong C$,
and $\mathcal{F}(\mathcal{G}(D)) \cong D$ $*$-isomorphically. 
One can think of $\mathcal{F}$
as a mapping between  $W^{*}$-covers
of $M$ and $N$.
There is a natural ordering of $W^*$-covers 
of a dual operator algebra.
If $(A, j)$ and $(A',j')$ are $W^*$-covers of $M$,
we then define $(A,j)$ $\leq (A',j')$
if and only if there is a $w^*$-continuous
$*$-homomorphism $\pi : A'\to A$ such
that $\pi \circ j' = j$. It is an easy exercise
(using that the range of $\pi$ is $w^*$-closed) to check that $\pi$ is surjective.

\begin{theorem}
The correspondence $C \mapsto \mathcal{F}(C)$
is bijective and order preserving. 
\end{theorem}

\begin{proof}
From the above discussion, the bijectivity is clear.
Suppose $\phi : C_{1} \to C_{2}$ is a $w^{*}$-continuous
quotient $*$-homomorphism between two $W^{*}$-algebras
generated by $M$, such that $\phi|_M$ = $Id_{M}$. Then
by Corollary 2.4 in \cite{BK1}
$\tilde{\phi}$ = 
$Id_{Y} \otimes \phi \otimes Id_{X} : Y \otimes^{\sigma h}_{M} C_{1} \otimes^{\sigma h}_{M} X \to Y \otimes^{\sigma h}_{M} C_{2} \otimes^{\sigma h}_{M} X $ 
is a $w^{*}$-continuous completely contractive
 map with $w^{*}$-dense range, which equals 
the identity when restricted to the copy of $N$. 
It is easy to check that $\tilde{\phi}$ is a 
homomorphism on the $w^*$-dense subset
$Y \otimes C_1 \otimes X$.
Therefore by $w^*$-density,  $\tilde{\phi}$ is a homomorphism.  
Hence by  Proposition A.5.8 in \cite{DBbook},
$\tilde{\phi}$ is a $*$-homomorphism and is onto. Hence, $\phi$
is order preserving.
\end{proof}

\begin{corollary}
If $\mathcal{L}^{w}$ is the linking dual operator
algebra for a $weak^{*}$ Morita equivalence of dual operator
algebras $M$ and $N$, and if $\mathcal{L'}$ is the 
corresponding linking $W^{*}$-algebra of the weak 
Morita equivalence of $W^{*}$-algebras $W^{*}_{max}(M)$
and $W^{*}_{max}(N)$, then $W^{*}_{max}(\mathcal{L}^{w}) = \mathcal{L'}$.
\end{corollary}

\begin{proof}
Suppose $W^{*}_{max}(M)$ is normally and 
faithfully represented on $B(H)$ for
some Hilbert space $H$. Then, by Lemma 1.1 in \cite{BK1},
$H$ is a normal universal Hilbert $M$-module.
 Also $M$ is weak$^{*}$ Morita equivalent 
to $\mathcal{L}^{w}$, via the
dual bimodule $M \oplus^{c} Y$ (see Corollary 4.1 in \cite{BK1}).
 By  Theorem 3.10 in \cite{BK1}, this induces a  normal 
representation of $\mathcal{L}^{w}$  on the
Hilbert space $(M \oplus_{c} Y) \otimes^{\sigma h}_{M} H^{c}$. 
By  Proposition 4.2 in \cite{BK1} we have that
\begin{center}
 $(M \oplus^{c} Y) \otimes^{\sigma h}_{M} H^{c}$ $\cong$ $ (H \oplus K)^{c}$
\end{center}
\noindent  unitarily, where $K$ = $Y \otimes^{\sigma h}_{M} H^c$
and $K$ is also a normal universal Hilbert $N$-module (see e.g. remark on page 6 in \cite{BSo}). 
As in the proof of Theorem 5.2 in
\cite {BK1}, $W^{*}_{max}(\mathcal{L}^{w})$ may
be taken to be the  $W^{*}$-algebra generated by 
$\mathcal{L}^{w}$ in $B(H \oplus K)$, which is $\mathcal{L'}$.
 \end{proof}

The above corollary should have a  variant
valid for arbitrary $W^*$-covers which we hope to include in
\cite{UK}. That is, if $\mathcal{L'}$
is the corresponding linking $W^*$-algebra
of the weak Morita equivalence of arbitrary $W^*$-covers
then $\mathcal{L'}$ is a $W^*$-cover of $\mathcal{L}^w$.

\begin{proposition} \label{L}
If $(M,N,X,Y)$ is a subcontext of a $W^{*}$-Morita
context $(\mathcal{E},\mathcal{F},W,Z)$, then
\begin{enumerate}
\item
$X$ and $Y$ generate $W$ and $Z$ respectively as 
left dual operator modules; i.e., $W$ is the smallest 
$w^{*}$-closed left $\mathcal{E}$-submodule of $W$ containing 
$X$. Similar assertions hold as right dual operator
modules, by symmetry.
\item
The linking algebra $\mathcal{L}$ of  $(M,N,X,Y)$ 
generates the
linking $W^{*}$-algebra $\mathcal{L'}$ of $(\mathcal{E},\mathcal{F},W,Z)$.
\item
If $M$ or $N$ is a $W^{*}$-algebra, then $(M,N,X,Y)$
= $(\mathcal{E}, \mathcal{F}, W,Z)$.
\end{enumerate}
\end{proposition}

\begin{proof}
Since the pairing $[\cdot,\cdot] : Y \times X \to N$ has $w^{*}$-dense
range, we can pick a net $e_{t}$ in $N$ which is a sum of terms of the form
$[y,x]$, for $y \in Y$, $x \in X$, such that $e_{t}  \buildrel w^{*} \over   \to 1_{N}$.
Hence $w e_{t} \buildrel w^{*} \over   \to w$ for all $w \in W$. Thus, sums of 
terms of the form $w[y,x]$, for $w \in W, x \in X, y \in Y$ are $w^{*}$-dense 
in $W$. However, $w[y,x] = (w,y) x$ $\in$ $\mathcal{E}X$ which shows that 
$\mathcal{E}X$ is $w^{*}$-dense in $W$. Thus, $X$ generates $W$ as a left
dual operator $\mathcal{E}$-module. Assertions (2) and (3) follow
from (1). For example, if $M$ is a $W^*$-algebra,
then clearly $X$ = $W$. Since $Y$ generates $Z$
as a right dual operator module, we have
$Z = \overline{Y\mathcal{E}}^{w*} = \overline{Y M}^{w*} = Y$. 
Since the ranges of the
natural pairings $ Z \times W \to \mathcal{F}$
and $Y \times X \to N$ are weak$^*$-dense, 
this implies that $\mathcal{F}$ = $N$.
\end{proof}

\begin{theorem}
If $(M,N,X,Y)$ is a weak$^{*}$ Morita context
which is a subcontext of a W$^{*}$-Morita context 
$(\mathcal{E},\mathcal{F},W,Z)$, 
then it is a dilatable subcontext. 
\end{theorem}

\begin{proof}
By Proposition \ref{L}, $X$ and $Y$ generate $W$ and $Z$, respectively,
as left dual operator modules. Hence we have a $w^{*}$-continuous
complete contraction $\mathcal{E} \otimes^{\sigma h}_{M} X \to W$
with $w^{*}$-dense range. On the other hand,
 \begin{center}
                        $W \cong W \otimes^{\sigma h}_{N} N \cong W \otimes^{\sigma h}_{N} Y \otimes^{\sigma h}_{M} X  \cong (W \otimes^{\sigma h}_{N} Y) \otimes^{\sigma h}_{M} X$\\
  \end{center}
  \noindent  completely isometrically and $w^{*}$-homeomorphically. 
  However, the pairing $(\cdot, \cdot) :  W \times Y \to \mathcal{E}$ 
  determines a $w^{*}$-continuous complete
  contraction $W \otimes^{\sigma h}_{M} Y \to \mathcal{E}$, and so we obtain a 
  $w^{*}$-continuous complete
  contraction $W \to \mathcal{E} \otimes ^{\sigma h}_{M} X$. 
  Recall from \cite{BK1} that $N$ has an `approximate identity' of the form
   $\sum_{i=1}^{n_t}  [y_i^t,x_i^t]$. 
   Under the above identifications,\\

\noindent $w \mapsto w \otimes_{N} 1_N$ $\mapsto w \otimes_{N}  w^*$-lim$_t $ $\sum_{i=1} ^{n_t}  y_i^t \otimes_{M} x_i^t$
 $\mapsto w^*$-lim$_t  \sum_{i=1} ^{n_t}  (w \otimes_{N} y_i^t ) \otimes_{M} x_i^t$\\

$\mapsto w^*$-lim$_t  \sum_{i=1} ^{n_t} (w, y_i^t ) x_i^t$
$\mapsto w^*$-lim$_t  \sum_{i=1}^{n_t} w [y_i^t,x_i^t]$
$= w$. \\

 \noindent Hence, the composition of these maps
  \begin{center}
                                         $\mathcal{E} \otimes^{\sigma h}_{M} X \to W \to \mathcal{E} \otimes^{\sigma h}_{M} X$\\
 \end{center}                                        
\noindent is the identity map,
from which it follows that $W \cong \mathcal{E} \otimes^{\sigma h}_{M} X$.
Similarly $Z$ is the dilation of $Y$.                                                               
\end{proof}

\begin{theorem} \label{mai}
If $(M,N,X,Y)$ is a left dilatable maximal subcontext
of a W$^{*}$-context, then $M$ and $N$ are $weak^{*}$
Morita equivalent dual operator algebras.
Indeed, it also follows that $(M,N,X,Y)$ is a 
weak$^{*}$ Morita context. Conversely, every $weak^{*}$Morita 
equivalence  of dual operator algebras occurs 
in this way. That is, every weak$^{*}$ Morita context 
is a left dilatable maximal subcontext of a W$^{*}$-Morita context.  
\end{theorem}

\begin{proof}
 Every weak$^{*}$ Morita context 
is a left dilatable maximal subcontext of a W$^{*}$-Morita context is proved in Theorem 5.2 in \cite{BK1}. For the converse,
let $\mathcal{C}$ and $\mathcal{D}$ be the usual
maximal $W^{*}$-algebras of $M$ and $N$ respectively,
and let $(M,N,X,Y)$ be a left dilatable  subcontext of $(\mathcal{C},\mathcal{D},W,Z)$.
 Using Lemmas \ref{C} and \ref{d}, we have
 \begin{center}
 $Y \otimes^{\sigma h}_{M} X \subset (\mathcal{D} \otimes^{\sigma h}_{N} Y) \otimes ^{\sigma h}_{M} X  
\cong Z \otimes ^{\sigma h}_{M} X \cong (Z \otimes^{\sigma h}_{\mathcal{C}}  \mathcal{C}) \otimes^{\sigma h}_{M} X \cong Z \otimes^{\sigma h}_{\mathcal{C}} W \cong \mathcal{D}$,
  \end{center}
 
 \noindent complete isometrically and $w^{*}$-homeomorphically. On the
 other hand, we have the canonical $w^{*}$-continuous complete contraction
 \begin{center}
 $Y \otimes^{\sigma h}_{M} X \to N \subset \mathcal{D}$
 \end{center}
 \noindent coming from the restricted pairing in  Definition \ref{PP} (2). 
It is easy to check that
 the composition of maps in these two sequences agree. 
Thus the canonical map
  $ Y \otimes^{\sigma h}_{M} X \to N$ is a 
$w^*$-continuous completely isometric isomorphism.
   Similarly,   
$ X \otimes^{\sigma h}_{N} Y \to M$ is a $w^*$-continuous
completely isometric isomorphism.
Hence by the Krein-Smulian theorem,
$ X \otimes^{\sigma h}_{N} Y $ $\cong$ $ M$
and $ Y \otimes^{\sigma h}_{M} X $ $\cong$  $N$ 
completely isometrically and $w^*$-homeomorphically.
Thus $M$ and $N$ 
are weak$^{*}$ Morita equivalent dual operator algebras.
\end{proof}

\section{The Main Theorem}

\begin{definition} \label{D}
Two dual operator algebras $M$ and $N$ are {\em (left)
dual operator Morita equivalent} if there exist completely 
contractive functors $F :$ $_{M} \mathcal{R} \to$ $_{N} \mathcal{R}
$  and  $G :$ $_{N} \mathcal{R} \to $ $_{M} \mathcal{R}$
which are weak$^{*}$-continuous on
morphism spaces, such that $FG$ $\cong$ $Id$ and $GF$ $\cong$ $Id$ completely isometrically.
Such $F$ and $G$ will be
called {\em dual operator equivalence functors}.
\end{definition}

Note that by Corollary 3.5.10 in \cite{DBbook},
$CB_{M}(V,W)$ for $V, W \in$ $_{M}\mathcal{R}$ is a dual operator space,
but $CB^{\sigma}_{M}(V,W)$ is not a $w^*$-closed subspace
of $CB_{M}(V,W)$. In the above definition, by the functor $F$
being $w^*$-continuous on morphism spaces, we mean that
if $(f_t) \subseteq   CB^{\sigma}_{M}(V,W)$, 
$f_t \buildrel w^* \over \to f$ in $CB_{M}(V,W)$,
and if $f$  also lies in $CB^{\sigma}_{M}(V,W)$, then $F(f_t) \buildrel w^* \over  \to F(f)$ in $CB_{N}(F(V), F(W))$. Similarly for the functor $G$. 
We also assume that the natural transformations
coming from $GF \cong Id$ and $FG \cong Id$ are weak$^*$-continuous
in the sense that for all $V \in$  $_{M} \mathcal{R}$,
the natural transformation $w_V : GF(V) \to V$
is a weak$^*$-continuous map. Similarly for $FG \cong Id$.

There is an obvious analogue to `right dual operator
Morita equivalence', where we are concerned with right 
dual operator modules. Throughout, we write
$\mathcal{C}$ and $\mathcal{D}$ for $W^*_{max}(M)$ 
and $W^*_{max}(N)$ respectively.

We now state our main theorem:

\begin{theorem} \label{main}
Two dual operator algebras are weak$^*$ Morita equivalent
if and only if they are left dual operator Morita equivalent,
and if and only if they are right dual operator
Morita equivalent.
Suppose that $F$ and $G$ are the left dual operator
equivalence functors, and set $Y$ = $F(M)$ and $X$ = $G(N)$.
Then $X$ is a weak$^*$ Morita equivalence $M$-$N$-bimodule. 
Similarly $Y$ is a weak$^*$ Morita equivalence
$N$-$M$-bimodule; that is, $(M, N, X, Y)$ is a weak$^*$ Morita context. 
Moreover, $F(V)$ $\cong$ $ Y \otimes^{\sigma h}_{M} V$
completely isometrically and weak$^{*}$-homeomorphically
(as dual operator $N$-modules) for all $V$ $\in$ $_{M} \mathcal{R}$.
Thus, $F$ $\cong$ $ Y \otimes^{\sigma h}_{M} -$ and $G$ $\cong$ $X \otimes^{\sigma h}_{N}-$ completely isometrically.
Also $F$ and $G$ restrict to equivalences of 
the subcategory $_{M} \mathcal{H}$ with
$_{N} \mathcal{H}$, the subcategory $_{\mathcal{C}}\mathcal{H}$ with $_{\mathcal{D}} \mathcal{H}$, and the subcategory $_{\mathcal{C}}\mathcal{R}$ with $_{\mathcal{D}} \mathcal{R}$.  
\end{theorem}

We will use techniques similar to those of \cite{DB3} and \cite{DB4}
to prove our main theorem. Mostly this involves the change of
tensor product and modification of arguments in the
present setting of weak$^*$-topology.

The following lemmas will be very useful to us. Their proofs
are almost identical to analogous results 
in \cite{DB3} and therefore are omitted.

\begin{lemma} \label{A}
Let $V$ $\in$ $_{M} \mathcal{R}$. Then $v$ $\mapsto$ $ r_v$ 
is a $w^{*}$-continuous complete isometry of $V$ onto $ CB_{M}(M,V)$.
In this case, $  CB_{M}(M,V)$ =  $   CB^{\sigma}  _{M}(M,V)$ i.e. $V$ $\cong$ $  
 CB^{\sigma} _{M}(M,V)$
completely isometrically and $w^*$-homeomorphically.
\end{lemma}

\begin{lemma}\label{B}
If $V$, $V'$ $\in$ $_{M}\mathcal{R}$ then the map
$T$ $\mapsto$ $F(T)$ gives a completely
isometric surjective linear isomorphism $  CB^{\sigma}_{M}(V,V')$ $\cong$
$ CB^{\sigma} _{N}(F(V),F(V'))$. If $V$ = $V'$, then this map
is a completely isometric surjective homomorphism. 
\end{lemma}

\begin{lemma} \label{qw}
For any $V$ $\in$ $_{M} \mathcal{R}$, we have $F(R_m(V))$ $\cong$
$R_m(F(V))$ and $F(C_m(V))$ $\cong$ $C_m(F(V))$ completely
isometrically.
\end{lemma}

\begin{lemma} \label{t}
The functors $F$ and $G$ restrict to a completely isometric
functorial equivalence of the subcategories  $_{M} \mathcal{H}$
and $_{N} \mathcal{H}$.
\end{lemma}

\begin{proof}

Let $H \in$ $ _{M} \mathcal{H}$. Recall that $H$ with its column Hilbert space structure 
$H^c$ is a left dual operator $M$-module. We need to show that $K = F(H^c) \in$ $ _{N}\mathcal{H}$
or equivalently $F(H^c)$ is a column Hilbert space. 
For any dual operator space $X$ and  $m \in \mathbb{N}$, we have $X \otimes_h C_m = X \otimes^{\sigma h} C_m$. Hence by Proposition 2.4 in \cite{DB3}, it suffices to show 
that the identity map $K \otimes_{min} C_m \to K \otimes^{\sigma h} C_m$
is a complete contraction for all $m \in \mathbb{N}$.
Since all operator space tensor products coincide for Hilbert column spaces, 
we have 
$C_m(H^c) \cong H^c \otimes_{min} C_m \cong H^c \otimes_h C_m \cong H^c \otimes^{\sigma h} C_m$. 
Thus \begin{eqnarray*}
K \otimes_{min} C_m &\cong&  C_m(F(H^c)) \\
&\cong &F(C_m(H^c)) \\
&\cong&  F(H^c \otimes^{\sigma h} C_m) \\
&\cong& F(G(K) \otimes^{\sigma h} C_m)  
\end{eqnarray*}
using Lemma \ref{qw} and $G(K) \cong H^c$.
Also, using Lemma \ref{A} and Lemma \ref{B} we have
\begin{eqnarray*}
G(K) &\cong& CB_{M}(M,G(K)) \\
&\cong& CB_{N}^{\sigma}(Y, FG(K) )\\
&\cong& CB_{N}^{\sigma}(Y, K).
\end{eqnarray*}
By Lemma 2.3 in \cite{BK1} we get a complete contraction
$G(K) \otimes^{\sigma h} C_m \to CB_{N}^{\sigma} (Y ,K) \otimes^{\sigma h} C_m$.
Now $CB_{N}^{\sigma} (Y ,K) \otimes^{\sigma h} C_m$
$\to CB^{\sigma}_{N}(Y, K \otimes^{\sigma h} C_m) 
: T \otimes z \mapsto y \mapsto T(y) \otimes z$
for $T \in CB_{N}^{\sigma} (Y ,K)$ and $z \in C_m$,
is a complete contraction.  
Again using Lemma \ref{A} and Lemma \ref{B},
we have $CB^{\sigma}_{N}(Y, K \otimes^{\sigma h} C_m) \cong CB^{\sigma}_{M} (M, G(K \otimes^{\sigma h} C_m)) \cong G(K \otimes^{\sigma h} C_m)$. Taking the composition of above maps
gives a complete contraction $G(K) \otimes^{\sigma h} C_m \to G(K \otimes^{\sigma h} C_m)$. 
Applying $F$ to this map, we get a complete contraction
 $F(G(K) \otimes^{\sigma h} C_m) \to K \otimes^{\sigma h} C_m$.
 This together with $K \otimes_{min} C_m \cong  F(G(K) \otimes^{\sigma h} C_m)  $
 gives the required complete contraction  $K \otimes_{min} C_m \to K \otimes^{\sigma h} C_m$.
 \end{proof}

\begin{corollary} \label{p}
The functors $F$ and $G$ restrict to a 
completely isometric equivalence of  $_{\mathcal{C}}\mathcal{H}$
and  $_{\mathcal{D}}\mathcal{H}$.
\end{corollary}

\begin{proof}
This is Corollary \ref{n} proved earlier. 
\end{proof}

Also, this restricted
equivalence is a normal $*$-equivalence in the sense of Rieffel
\cite{R}, and so $\mathcal{C}$ and $\mathcal{D}$ are weak Morita
equivalent in the sense of  Definition 7.4 in \cite{R}.

\begin{lemma} \label{U}
For a dual operator $M$-module $V$,
the canonical map $ \tau_{V} : Y \otimes V \to F(V)$
given by $y \otimes v \mapsto F(r_v)(y)$ is separately $w^*$-continuous 
and extends uniquely to a completely contractive
map on $Y \otimes^{\sigma h}_{M} V$. Moreover, this map has $w^*$-dense
range.
\end{lemma}

\begin{proof}
Since the
functor $F$ is $w^*$-continuous on morphism spaces,
it is easy to check that $\tau _V : Y \times V \to F(V)$
is a separately $w^*$-continuous bilinear map.
To see that $\tau_V$ has $w^*$-dense range, suppose the contrary.
Let $Z$ = $F(V)/N$ where $N$ = $\overline{Range(\tau_V)}^{w*}$
and let $Q : F(V) \to Z$
be the nonzero $w^*$-continuous quotient map. 
Then $G(Q): G(F(V)) \to G(Z)$ 
is nonzero. Thus there exists $v$ $\in$ $V$ such that
$G(Q)w_V^{-1}r_v$ $\neq$ $0$ as a map on
$M$, where $w_V$ is the $w^*$-continuous completely isometric 
natural transformation
 $GF(V) \to V$ coming from $GF \cong Id$.
Hence $FG(Q)F(w_V^{-1})F(r_v)$ $\neq$ $0$, and thus $QTF(r_v)$
$\neq$ $0$ for some $w^*$-continuous module map 
$T : F(V) \to F(V)$ since $w_V^{-1}$
is $w^*$-continuous by the Krein-Smulian theorem. By
Lemma \ref{B}, $T$ = $F(S)$ for some $w^*$-continuous
module map $S : V \to V$, 
so that $QF(r_{v'})$ $\neq$ $0$ for $v' = S(v)$ $\in$ $V$.
Hence $Q \circ \tau_V$ $\neq$ $0$, which is a contradiction. 
Again as in the proof of 
Lemma 2.6 in \cite{DB4}, $\tau_V$ is a complete contraction. 
Thus, $\tau_V$
is a separately $w^*$-continuous completely
contractive bilinear map. The result follows from the 
universal property of $Y \otimes^{\sigma h}_{M} V$.
\end{proof}

Let $(M, N, \mathcal{C}, \mathcal{D}, F, G, X, Y)$ be as above.
We let $H$ $\in$ $_{M} \mathcal{H}$
 be the Hilbert space of the normal universal representation of $\mathcal{C}$
and let $K$ = $F(H)$. By Lemma \ref{t} and Corollary \ref{p},
$F$ and $G$ restrict to equivalences of $_{M} \mathcal{H}$ with
$_{N} \mathcal{H}$, and restrict further to  normal
$*$-equivalences of   $_{\mathcal{C}} \mathcal{H}$ with
 $_{\mathcal{D}}\mathcal{H}$. 
 By Proposition 1.3 in \cite{R} , $\mathcal{D}$
 acts faithfully on $K$.
 Hence, we can regard $\mathcal{D}$ as a 
 subalgebra of $B(K)$.
 Define $Z$ = $F(\mathcal{C})$
and $W$ = $G(\mathcal{D})$.

From Lemma \ref{U}, with $V$ = $M$, it follows that
$Y$ is a right dual operator $M$-module with module
action $y \cdot m$ = $F(r_{m})(y)$, for $y \in Y$, $m \in M$ 
and $r_{m} : M \to M$ :  $c \mapsto cm$ is simply
right multiplication by $m$. Similarly,
$X$ is a right dual operator $N$-module, 
and $Z$ and $W$ are dual operator 
$N$-$\mathcal{C}$- and $M$-$\mathcal{D}$-bimodules
respectively. The inclusion $i$ of
$M$ in $\mathcal{C}$ induces a completely
contractive $w^*$-continuous inclusion
$F(i)$ of $Y$ in $Z$. One can check
that $F(i)$ is a $N$-$M$-module map. By Lemma ~ \ref{T} below
and its proof, it is easy to see that  $F(i)$ is a complete isometry.
Hence we may regard $Y$ as a $w^*$-closed $N$-$M$-submodule
of $Z$ and similarly $X$ may be regarded as a $w^*$-closed $M$-$N$-submodule of $W$.

With $V$ = $X$ in Lemma \ref{U}, there is a 
left $N$-module map $ Y \otimes X$ $\to$ $F(X)$
defined by $y \otimes x \mapsto F(r_{x})(y)$. Since
$F(X)$ = $FG(N)$ $\cong$ $N$, we get a left $N$-module map
$[.] : Y \otimes X \to N$. In a similar
way we get a module map $(.) : X \otimes Y \to M$.
In what follows we may use the same notation for the
unlinearized bilinear maps, so for example we may use the symbol
$[y, x]$ for $[y \otimes x]$. 
These maps $(.)$ and $[.]$ have natural extensions
to $Y \otimes W \to \mathcal{D}$ and 
$X \otimes Z \to \mathcal{C}$ respectively,
which we denote by the same symbols.
Namely, $[y,w]$ is defined via $\tau_{W}$ for $y \in Y$ and $w \in W$. 
By Lemma \ref{U}, these maps
have weak$^{*}$-dense ranges.  

\begin{lemma} \label{T}
The canonical maps $X \to  CB^{\sigma} _{N}(Y,N)$ and $Y \to CB^{\sigma} _{M}(X,M)$,
induced by $[.]$ and $(.)$  respectively, 
are completely isometrically isomorphic. 
Similarly, the extended maps $W \to  CB^{\sigma} _{N}(Y,\mathcal{D})$
and $Z \to  CB^{\sigma}_{M}(X,\mathcal{C})$ are complete isometries.
\end{lemma}

\begin{proof} 
By Lemma \ref{A} and Lemma \ref{B},
we have $X$ $\cong$ $CB^{\sigma} _{M}(M,X)$ $\cong$ $CB^{\sigma} _{N}(Y,F(X))$
$\cong$ $CB^{\sigma} _{N}(Y,N)$ completely isometrically. Taking the composition
of these maps  shows that $x \in X$ corresponds to the map
$y \mapsto [y,x]$ in $CB^{\sigma} _{N}(Y,N)$. Similarly for the other maps. 
 \end{proof}

Next consider maps $\phi : Z \to B(H,K)$, and $\rho: W \to B(K,H)$
defined by $ \phi(z)(\zeta)$ = $F(r_{\zeta})(z)$, and $\rho(w)(\eta)$
= $\omega_{H} G(r_{\eta})(w)$, for $\zeta \in H$ and $\eta \in K$
where $\omega_{H} : GF(H) \to H$ is the $w^*$-continuous
$M$-module map coming from the natural transformation $GF$ $\cong$ $Id$.
Again $r_{\zeta} : \mathcal{C} \to H$ and $r_{\eta} : \mathcal{D}
 \to K$ are the 
obvious right multiplications. As $\omega_{H}$ is an isometric onto map
between Hilbert spaces, $\omega_{H}$ is unitary and hence 
also a $\mathcal{C}$-module map by Corollary \ref{cor}.
 One can check that:
\begin{equation} \label{*}
\rho(x) \phi(z) = (x, z) \ \      and    \ \    \phi(y) \rho(w) = [y,w]V 
\end{equation}

\noindent for all $x$ $\in$ $X$, $y$ $\in$ $Y$, $z \in Z$, $w \in W$ and where $V$ $\in$ $B(K)$ 
is a unitary operator in $\mathcal{D} '$ composed of two natural transformations.
A similar calculation as in Lemma ~ 4.3 in \cite{DB4}, shows that 
the unitary $V$ is in the center of
 $\mathcal{D}$, hence  $\phi(y) \rho(w)$ $\in \mathcal{D}$ for all $y \in Y$ and $w \in W$.

\begin{lemma} \label{x}
The map $\phi$ (respectively $\rho$) is a completely isometric
$w^*$-continuous $N$-$\mathcal{C}$-module map 
(respectively $M$-$\mathcal{D}$-module map). Moreover,
$\phi(z_1)^* \phi(z_2)$ $\in$ $\mathcal{C}$ for all $z_1, z_2$ $\in$
$Z$, and $\rho(w_1)^*\rho(w_2)$ $\in$ $\mathcal{D}$, for all $w_1, w_2$ $\in$ $W$. 
\end{lemma}

\begin{proof}
We will prove that the
maps $\phi$ and $\rho$ are $w^*$-continuous.
The rest of the assertions follow as in  Lemma 4.2 in \cite{DB4} 
and by von Neumann's double commutant theorem. To see that $\phi$ 
is $w^*$-continuous, let $(z_t)$ be a bounded net in $Z$
such that $z_{t}$  $ \buildrel w^* \over \to  z$ in $Z$. 
For $\zeta \in H$,  we have $F(r_{\zeta})$ $\in$ $ CB^{\sigma} _{N}(Z,K)$.
Hence $F(r_{\zeta})(z_{t})$ $\to$ $F(r_{\zeta})(z)$ weakly. That is,
$\phi(z_{t})$ $\to$ $\phi(z)$ in the WOT and it follows that
$\phi$ is weak$^*$-continuous. A similar argument works for $\rho$. 
\end{proof}

We will follow the approach of \cite{DB3} to prove the selfadjoint
analogue of our main theorem, which  involves
a change of the tensor product. Nonetheless, for completeness
we will give the proof.

\begin{theorem} \label{sec}
Two $W^*$-algebras $A$ and $B$ are weakly Morita equivalent
in the sense of Rieffel 
if and only if they are dual operator Morita equivalent 
in the sense of Definition \ref{D}. Suppose that $F$ and $G$
are the dual operator equivalence functors, and set $Z$ = $F(A)$
and  $W$ = $G(B)$. Then, $W$ is a $W^*$-equivalence 
$A$-$B$-bimodule, $Z$ is a $W^*$-equivalence 
$B$-$A$-bimodule, and $Z$ is unitarily and
$w^*$-homeomorphically isomorphic
to the conjugate $W^*$-bimodule $\overline{W}$ of $W$.
 Moreover, $F(V)$ $\cong$ $ Z \otimes^{\sigma h}_{A} V$
completely isometrically and weak$^{*}$-homeomorphically
(as dual operator $B$-modules) for all $V$ $\in$ $_{A} \mathcal{R}$.
Thus $F$ $\cong$ $ Z \otimes^{\sigma h}_{A} -$ and $G$ $\cong$ $W \otimes^{\sigma h}_{B}-$ completely isometrically.
Also $F$ and $G$ restrict to equivalences of 
the subcategory $_{A} \mathcal{H}$ with
$_{B} \mathcal{H}$.
\end{theorem}

\begin{proof}
In \cite{BK1} we saw that the weakly Morita equivalent
$W^*$-algebras (in the sense of Rieffel) are
weak$^*$ Morita equivalent. Hence by Theorem 3.5 in \cite{BK1},
they have equivalent categories of dual operator modules and the assertion
about the form of the functors also holds.

For the other direction, observe that by Corollary \ref{t},
the functors $F$ and $G$ restrict to a completely isometric
equivalence of $_{A} \mathcal{H}$ and $_{B} \mathcal{H}$. Hence, by  
Definition 7.4 in \cite{R}, $A$ and $B$ are weakly Morita 
equivalent in the sense of Rieffel.
 We will follow \cite{DB3} to prove 
rest of the assertions.


By the polarization identity and Lemma \ref{x}, 
$W$ is a right $C^*$-module over $B$ with inner product
$\langle w_1,  w_2 \rangle_{B}$ = $\rho(w_1)^{*} \rho(w_2)$, 
for $w_1, w_2 \in W$. Similarly,
$W$ is a left $C^*$-module over $A$ by setting 
$_{A} \langle w_1,w_2 \rangle$ = $\rho(w_1) \rho(w_2)^*$.
To see that this inner product lies in $A$, note that, since the
range of $(.)$ is $w^*$-dense in $A$,  we can
choose a net in $A$ of the form $e_{\alpha}$ = $\sum_{k=1}^{n(\alpha)} (w_k, z_k)$
= $\sum_{k=1} ^{n(\alpha)} \rho(w_k) \phi(z_k)  $ 
where $z_k \in Z$ and $w_k \in W$,
such that $e_{\alpha}$  $ \buildrel w^* \over \to  1_{A}$. 
Then, $e_{\alpha}^*$  $ \buildrel w^* \over \to  1_{A}$.
 Since $\rho$ is a weak$^*$-continuous
 $A$-module map, $\rho(w)^*$ = $w^*$-lim$_{\alpha}$ $\rho(e_{\alpha}^* w)^*$
= $w^*$-lim$_{\alpha}$ $ \rho(w)^* e_{\alpha}$, it follows that
$\rho(w) \rho(w)^*$ is a weak$^*$ limit
of finite sums of terms of the form 
$\rho(w)( \rho(w)^* \rho(w_k)) \phi(z_k)$ = $\rho(w) \phi(b z_k) $
= $(w, b z_k)$ $\in$ $A$, where $b$ = $ \rho(w)^* \rho(w_k)$ $\in$ $B$.
Thus $\rho(w) \rho(w)^*$ $\in$ $A$.
By the polarization identity $\rho(w_1) \rho(w_2)^*$ $\in$ $A$. 
Similarly, $Z$ is both a left and a right $C^*$-module. To see that $Z$
is a $w^*$-full right $C^*$-module over $A$, 
 rechoose a net
in $A$ of the form  $e_{\alpha}$ = $\sum_{k=1} ^{n(\alpha)} \rho(w_k)\phi(z_k)$
such that $e_{\alpha}$ $\to$ $I_H$ strongly, 
so that $e_{\alpha}^*e_{\alpha}$ $\to$ $I_{H}$ weak$^*$ as done  in Theorem 4.4
in \cite{BK1}.
However $e_{\alpha}^*e_{\alpha}$ =  $ \sum_{k,l} \phi(z_k)^*  b_{kl} \phi(z_l)$
where $b_{kl}$ = $ \rho(w_k)^* \rho(w_l)$ $\in$ $B$.
Since $P$ = $[b_{kl}]$ is a positive matrix, it has a square root
$R$ = $[r_{ij}]$, with $r_{ij}$ $\in$ $B$. Thus 
$e_{\alpha}^* e_{\alpha}$ = $ \sum_{k} \phi(z_k^{\alpha})^* \phi(z_k^{\alpha})$
where $z_{k}^{\alpha}$ = $\sum_{j} r_{kj} z_{j}$.
From this one can easily deduce that 
the $A$-valued inner product on $Z$
has $w^*$-dense range. Similarly $Z$ is a weak$^*$-full left $C^*$-module
over $B$. Similarly for $W$. Since $\rho$ and $\phi$ are
$w^*$-continuous, the inner products are separately $w^*$-continuous.
Hence, by Lemma 8.5.4 in \cite{DBbook}, $W$ and $Z$ are
$W^*$-equivalence bimodules, implementing
the weak Morita equivalence of $A$ and $B$. Note that 
by Corollary 8.5.8
in \cite{DBbook}, 
$CB_{A}(W,A)$ = $ CB^{\sigma} _{A}(W,A)$. Thus by (8.18) in
\cite{DBbook} and Lemma \ref{T}, 
$Z$ $\cong$ $\overline{W}$ completely isometrically. 

Let $V \in\  _{\mathcal{A}} \mathcal{R}$. By Lemma \ref{A},  Lemma \ref{B}
above, Theorem 2.8 in \cite{BK2}, and the
fact that $Z$ $\cong$ $\overline{W}$, we have the following sequence of isomorphisms:
\begin{center}
$F(V) \cong \  CB^{\sigma} _{B}(B,F(V))  \cong  CB^{\sigma} _{A}(W,V)
 \cong Z \otimes^{\sigma h}_A V$
\end{center}
as left dual operator $B$-modules. Thus the conclusions of the theorem all hold.
\end{proof}
\medskip

Now we will come back to the setting
where $M$ and $N$ are dual operator algebras
and $\mathcal{C}$ and $\mathcal{D}$ are maximal
$W^*$-algebras generated by $M$ and $N$ respectively.

\begin{theorem} \label{up}
The $W^*$-algebras $\mathcal{C}$ and $\mathcal{D}$
are weakly Morita equivalent. In fact $Z$, which is a dual operator 
$N$-$\mathcal{C}$-bimodule, is a 
$W^*$-equivalence $\mathcal{D}$-$\mathcal{C}$-bimodule.
Similarly, $W$ is a $W^*$-equivalence $\mathcal{C}$-$\mathcal{D}$-bimodule, 
  and 
$W$ is unitarily 
and $w^*$-homeomorphically isomorphic
to the conjugate $W^*$-bimodule $\overline{Z}$ of $Z$ (and as dual operator bimodules).
\end{theorem}

\begin{proof}
By  Lemma \ref{x}, it follows that $\rho(W)$ is a $w^*$-closed TRO
(a closed subspace $Z \subset B(K,H)$ with $ZZ^*Z \subset Z$).
Hence, by 8.5.11 in \cite{DBbook} and 
Lemma \ref{x},  $W$ (or equivalently $\rho(W)) $
is a right $W^*$-module over $\mathcal{D}$
with inner product  
$\langle w_1,w_2 \rangle _{\mathcal{D}}$ = $\rho(w_1)^* \rho(w_2)$.  
Since $\rho$ is a complete isometry,
the induced norm on $W$ coming from the inner
product coincides with the usual norm.
 Similarly $Z$ is a right
$W^*$-module over $\mathcal{C}$. 
Also, $W$ (or equivalently $\rho(W)$) is a $w^*$-full left 
$W^*$-module over $\mathcal{E}$ =  weak$^*$ closure of $\rho(W) \rho(W)^*$,
with the obvious inner product $_{\mathcal{E}} \langle w_1, w_2 \rangle$ = $ \rho(w_1) \rho(w_2)^*$.
We will show that $\mathcal{E}$ = $\mathcal{C}$. Analogous statements
hold for $\mathcal{D}$ and $\phi$. It will be understood
that whatever a property is proved for $W$, by symmetry, the matching
assertions for $Z$ hold.

Let $\mathcal{L}^w$ be the linking $W^*$-algebra for the
right $W^*$-module $W$, viewed as a weak$^*$-closed
subalgebra of $B(H \oplus K)$. We let $\mathcal{A}$ = weak$^*$ closure of $\rho(W) \phi(Y)$.
It is easy to check, using the fact that
$\phi(Y) \rho(W) \in \mathcal{D}$ (see above Lemma \ref{x})
and Lemma ~ \ref{x},
that $\mathcal{A}$ is a dual operator
algebra.
By the last assertion of Lemma \ref{U} and \eqref{*}, 
$M = \overline{\rho(X) \phi(Y)}^{w*} \subseteq \mathcal{A}$
and the identity of $M$ is an identity of $\mathcal{A}$. 
We let $\mathcal{U}$  be the weak$^*$ closure
of $\mathcal{D} \phi(Y)$, and we define $\mathcal{L}$ to be the
following subset of $B(H \oplus K)$:
$$\left[
\begin{array}{ccl}
\mathcal{A} &  \rho(W) \\
\mathcal{U} & \mathcal{D}
\end{array}
\right] $$
Using \eqref{*} and Lemma \ref{x}, it is easy to check that $\mathcal{L}$ is a 
subalgebra of $B(H \oplus K)$. By explicit computation and
Cohen's factorization theorem,  
$\mathcal{L}^w \mathcal{L}$ = $\mathcal{L}$ and 
$\mathcal{L} \mathcal{L}^w$ = $\mathcal{L}^w$. 
Indeed, by Lemma \ref{x} and the fact that $\rho(W)$ is a TRO,
it follows that $\mathcal{L}^w \mathcal{L} \subseteq \mathcal{L}$.
Again by using \eqref{*}, Lemma \ref{x} and the fact that $\rho(W)^*$
is a left $W^*$-module over $\mathcal{D}$, it follows that 
$\mathcal{L} \mathcal{L}^w \subseteq \mathcal{L}^w$.
As $\rho(W)$ is a right $W^*$-module over $\mathcal{D}$
so $\rho(W)$ is a nondegenerate $D$-module (see e.g. 8.1.3 in \cite{DBbook}),
hence $\rho(W) = \rho(W) \mathcal{D}$ by Cohen's factorization theorem (A.6.2 in \cite{DBbook}).
For the same reason, $\rho(W) = \rho(W) \rho(W)^* \rho(W)$.
Now one can easily check that $\mathcal{L} \subseteq \mathcal{L}^w \mathcal{L}$
and similarly $\mathcal{L}^w \subseteq \mathcal{L} \mathcal{L}^w$.
Hence $\mathcal{L}^w \mathcal{L}$ = $\mathcal{L}$ and 
$\mathcal{L} \mathcal{L}^w$ = $\mathcal{L}^w$. 
Therefore,
we conclude that 
$\mathcal{L}^w$ = $\mathcal{L}$. Comparing corners of these
algebras gives $\mathcal{E}$ = $\mathcal{A}$ and $\mathcal{U}$ =
$ \rho(W)^*$. Thus, $M$ $\subseteq$ $\mathcal{E}$, from which it
follows that $\mathcal{C}$ $\subseteq$ $\mathcal{E}$, since
$\mathcal{C}$ is the $W^*$-algebra generated by $M$ in $B(H)$.
Thus $\rho(W)$ is a left $\mathcal{C}$-module,
so $W$ can be made into a left $\mathcal{C}$-module 
in a unique way (by Theorem \ref{two*}).
Also by Corollary
\ref{cor}, $\rho$ is a left
$\mathcal{C}$-module map. By symmetry, $Z$ is a left $\mathcal{D}$-module
and $\phi$ is a $\mathcal{D}$-module map, so that 
\noindent $\rho(W)^*$ = $\mathcal{U}$ = $\overline{\mathcal{D} \phi(Y)}^{w^*}$ $\subset$ $\phi(Z)$.
By symmetry,
$\phi(Z)^*$ $\subset$ $\rho(W)$, so that $\rho(W)^*$ = $\phi(Z)$.
Since, $\phi(Z) = \overline{ \mathcal{D} \phi(Y)}^{w*}$, by symmetry,
$\rho(W) = \overline{ \mathcal{C}  \rho(X) }^{w*}$.
Also, $ \rho(W) \phi(Y)$ $\subset$   $\overline{ \mathcal{C}  \rho(X)\phi(Y)}^{w^*}$ $\subset$ $\mathcal{C}$
and thus  $\mathcal{E}$ = $\mathcal{A}$ $\subset$ $\mathcal{C}$.
Thus $\mathcal{E}$ = $\mathcal{A}$ = $\mathcal{C}$,   
and that $\mathcal{D}$ = $\overline{\phi(Z) \phi(Z)^*}^{w^*}$ =  $\overline{\rho(W)^* \rho(W)}^{w^*}$.
This proves the theorem. 
\end{proof}

\section{$W^*-$restrictable equivalences}

\begin{definition}
We say that a dual operator equivalence
functor $F$ is {\em $W^*$-restrictable},
if $F$ restricts to a functor from
$_{\mathcal{C}} \mathcal{R}$ into 
$_{\mathcal{D}} \mathcal{R}$. 
\end{definition}
We prove our main theorem under the assumption
that the functors $F$ and $G$ are $W^*$-restrictable. Later
we will prove that this condition is automatic; i.e., the functors
$F$ and $G$ are automatically $W^*$-restrictable.\\

\begin{remark}
 The canonical equivalence
functors 
coming from a given weak$^*$ Morita 
equivalence are $W^*$-restrictable.
Suppose that $M$ and $N$ are weak$^*$ Morita 
equivalent and let $(M, N, X, Y)$ be a weak$^*$ Morita  context. Then from Theorem ~ 5.2 in \cite{BK1}
we know that $\mathcal{C}$ and $\mathcal{D}$
are weakly Morita equivalent $W^*$-algebras, with 
 $W^*$-equivalence  $\mathcal{D}$-$\mathcal{C}$-bimodule $Z = Y \otimes^{\sigma h}_{M} \mathcal{C}$.
From Theorem 3.5 in \cite{BK1},
$F(V) = Y \otimes^{\sigma h}_{M} V$, for
$V$ a dual operator  $M$-module.
However, 
if $V$ is a dual operator
$\mathcal{C}$-module,
$Y \otimes^{\sigma h}_{M} V$ $\cong$  
$Y \otimes^{\sigma h}_{M} \mathcal{C} \otimes^{\sigma h}_{\mathcal{C}} V$ $\cong$  
$Z \otimes^{\sigma h}_{\mathcal{C}} V$. Hence, $F$ restricted to 
$_{\mathcal{C}} \mathcal{R}$ is equivalent to $Z \otimes^{\sigma h}_{\mathcal{C}}-$,
and thus is $W^*$-restrictable.  
\end{remark}
\medskip

\begin{theorem} \label{M}
Suppose that the dual operator equivalence functors $F$ and $G$  
are $W^*$-restrictable. Then the conclusions of the Theorem \ref{main} all hold. 
\end{theorem}

\begin{proof}
Clearly, $F$ and $G$ gives a dual operator Morita equivalence of
$_{\mathcal{C}}\mathcal{R}$ and $_{\mathcal{D}}\mathcal{R}$  
when restricted to these subcategories. Set $Y$ = $F(M)$, $Z$ = $F(\mathcal{C})$,
$X$ = $G(N)$, and $W$ = $G(\mathcal{D})$ as before. 
By Theorem \ref{sec},  $\mathcal{C}$ and $\mathcal{D}$ are weakly Morita equivalent
von Neumann algebras with $Z$ and $W$ as $W^*$-equivalence
bimodules. From the discussion above Lemma \ref{T},  
$Y$ is a right dual operator $M$-module
and $X$ is a right dual operator $N$-module.
Also $Y$ is a $w^*$-closed $N$-$M$-submodule
of $Z$ and $X$ is a $w^*$-closed $M$-$N$-submodule of $W$.


For any left dual operator $\mathcal{C}$-module $X'$, we have the following
sequence of canonical complete isometries by Lemma \ref{A} and Lemma \ref{B}:
\begin{eqnarray*} 
 CB^{\sigma} _{M}(X,X') &\cong&  
CB^{\sigma}_{N}(N,F(X')) \\
 &\cong& F(X') \\
 &\cong&   CB^{\sigma}_{\mathcal{D}}(\mathcal{D}, F(X'))\\
&\cong& CB^{\sigma}_{\mathcal{C}}(W, X').
\end{eqnarray*} 
Hence, by  the discussion following Definition \ref{dilation},
and by  Lemma \ref{d}, we have
$W$ $\cong$ $\mathcal{C} \otimes^{\sigma h}_{M} X$
completely isometrically and as $\mathcal{C}$-modules. It
can be checked that this isometry is a right $N$-module map.
Similarly, $Z$ $\cong$ $\mathcal{D} \otimes^{\sigma h}_{N} Y$.

For any dual operator $M$-module $V$, we have, $Y \otimes^{\sigma h}_{M} V$ $\subset$
$(\mathcal{D} \otimes^{\sigma h}_{N}  Y ) \otimes^{\sigma h}_{M} V$  $\cong$ $Z \otimes^{\sigma h}_{M} V$ 
completely isometrically, since any dual operator module is
contained in its maximal dilation. On the other hand, using Lemma \ref{U}, 
Lemma \ref{B}, and Theorem \ref{sec}, respectively, we have the
following sequence of canonical completely contractive $N$-module maps:
\begin{center} 
$Y \otimes^{\sigma h}_{M} V \to F(V) 
 \to F(\mathcal{C} \otimes^{\sigma h}_{M} V)$
 $\cong$ $Z \otimes^{\sigma h}_{\mathcal{C}} (\mathcal{C} \otimes^{\sigma h}_{M} V)$
 $\cong$ $Z \otimes^{\sigma h}_{\mathcal{C}} V$. \\
 \end{center}
 
 The composition of the maps in this sequence coincides
 with the the composition of complete isometries in the last sequence.
 Hence, the canonical map $Y\otimes^{\sigma h}_{M} V$ $\to$ $F(V)$
 is a $w^*$-continuous complete isometry. Since this map has $w^*$-dense
 range, by the Krein-Smulian theorem, it is a complete isometric isomorphism.
 Thus $F(V)$ $\cong$ $Y \otimes^{\sigma h}_{M} V$, 
 and similarly $G(U)$ $\cong$ $X \otimes^{\sigma h}_{N} U$.
 Finally,  $M \cong GF(M) \cong X \otimes^{\sigma h}_{N} Y$, using Lemma 2.10 in \cite{BK1}
 and similarly $N \cong Y \otimes^{\sigma h}_{M} X$ completely isometrically
 and $w^*$-homeomorphically. 
\end{proof}

\begin{corollary}
Dual operator equivalence functors are
automatically $W^*$-restrictable.
\end{corollary}

\begin{proof}
Firstly, we will show that $W$ is the maximal dilation
of $X$, and $Z$ is the maximal dilation of $Y$. 
In  Theorem \ref{up}, we saw that
the set $\mathcal{U}$ equals $Z$.
This implies that $Y$ generates $Z$ as a left dual operator
$\mathcal{D}$-module. Similarly, $X$ generates $W$ as a 
left dual operator $\mathcal{C}$-module.

By Lemma \ref{A} and Lemma \ref{B}, we have the following sequence of maps
\begin{center}
$ CB^{\sigma} _{M}(X,H)$ $\cong$ $  CB^{\sigma} _{N}(N,K)$ $\cong$ $K$ $\cong$ 
$ CB^{\sigma} _{\mathcal{D}}(\mathcal{D},K) \to$ $CB^{\sigma}_{M}(W,H)$.
\end{center}
One can check that $\eta \in K$ corresponds under the last two maps in the sequence to the map $w \mapsto \rho(w)(\eta)$, which lies in 
$ CB^{\sigma} _{\mathcal{C}}(W,H)$, 
since $\rho$ is a left $\mathcal{C}$-module map. 
Thus, the composition $R$ of  the maps
in the above sequence has range contained in $ CB^{\sigma}_{\mathcal{C}}(W,H)$.
Also, $R$ is an inverse to the restriction map
$CB^{\sigma}_{\mathcal{C}}(W,H) \to $ $CB^{\sigma}_{M}(X,H)$.
Thus $CB^{\sigma}_{\mathcal{C}}(W,H)$ $\cong$  $ CB^{\sigma}_{M}(X,H)$.
Since $H$ is a normal universal representation
of  ${\mathcal{C}}$ (see the paragraph below Lemma ~ \ref{U}), it follows from
Theorem ~ \ref{wow}, that $W$ is the maximal dilation of $X$. Similarly
$Z$ is the maximal dilation of $Y$.   

Let $V \in \ _{\mathcal{C}} \mathcal{R}$. By Lemma \ref{A},  Lemma \ref{B},
 Definition \ref{dilation}, Theorem 2.8 in \cite{BK2}, and Theorem \ref{up},
 we have the following sequence of isomorphisms\\

$F(V) \cong   CB^{\sigma} _{N}(N,F(V)) \cong  CB^{\sigma}_{M}(X,V)
\cong  CB^{\sigma}_{\mathcal{C}}(W,V) \cong Z \otimes^{\sigma h}_{\mathcal{C}} V$,\\

\noindent as left dual operator $N$-modules. Since 
$Z \otimes^{\sigma h}_{\mathcal{C}} V$ is a left dual operator
$\mathcal{D}$-module, we see that $F(V)$ is a left dual operator 
$\mathcal{D}$-module and by Theorem \ref{two*}, this $\mathcal{D}$-module
action is unique.
Also by Corollary \ref{cor} the map 
$Z \otimes^{\sigma h}_{\mathcal{C}} V $ $\to$ $F(V)$ coming from 
the composition of the above isomorphisms is a $\mathcal{D}$-module map.
This map $Z \otimes^{\sigma h}_{\mathcal{C}} V $ $\to$ $F(V)$ is defined
analogously to the map $\tau_{V}$ defined in Lemma \ref{U}.
One can check that if $T: V_1 \to V_2$ is a morphism
in $_{\mathcal{C}} \mathcal{R}$, then the following diagram commutes:

$$\xymatrix{ Z \otimes^{\sigma h}_{\mathcal{C}} V_1 \ar[r] \ar[d]^{I_Z \otimes T}   &  F(V_1) \ar[d]^{F(T)}  \\         
   Z \otimes^{\sigma h}_{\mathcal{C}} V_2 \ar[r]  &   F(V_2) } $$

By Corollary 2.4 in \cite{BK1},
$I_Z \otimes T$ is a $w^*$-continuous $\mathcal{D}$-module map
and both the horizontal arrows above are  $w^*$-continuous
$\mathcal{D}$-module maps.
Hence, $F(T)$ is a $w^{*}$-continuous $\mathcal{D}$-module map; 
that is, $F(T)$ is a morphism in $_{\mathcal{D}} \mathcal{R}$.
Thus $F$ is $W^*$-restrictable. By  Theorem \ref{M}, our main
theorem is proved. 
\end{proof}


\subsection*{Acknowledgments} The author would like to thank David Blecher for 
many helpful comments, corrections, and suggestions during the preparation of this work.

\end{document}